\documentclass[letterpaper,11pt]{article}

\usepackage[margin=1in]{geometry}
\usepackage{amsmath,amssymb,amsfonts,amsthm,mathtools,dsfont}
\usepackage{graphicx,xcolor,booktabs,multirow,float,appendix,enumitem}
\usepackage{microtype}
\usepackage{url}
\usepackage[colorlinks=true,
            linkcolor=blue,
            citecolor=red,
            urlcolor=blue,
            anchorcolor=red,
            bookmarks=true,
            plainpages=false]{hyperref}
\usepackage[capitalize,nameinlink]{cleveref}
\usepackage{parskip}
\allowdisplaybreaks

\theoremstyle{plain}
\newtheorem{theorem}{Theorem}
\newtheorem{lemma}{Lemma}
\newtheorem{proposition}{Proposition}

\theoremstyle{definition}
\newtheorem{definition}{Definition}

\newtheorem{remark}{Remark}

\newcommand{\R}{\mathbb{R}}

\newcommand{\eps}{\varepsilon}
\newcommand{\dd}{\mathrm{d}}
\newcommand{\ci}{\perp\!\!\!\perp}

\renewcommand{\P}{\mathsf{P}}
\newcommand{\Q}{\mathsf{Q}}
\newcommand{\Expect}{\mathbb{E}}
\newcommand{\Prob}{\mathsf{P}}
\newcommand{\TV}{\operatorname{TV}}
\newcommand{\ind}[1]{\mathbf{1}_{\left\{#1\right\}}}

\DeclareMathOperator{\Unif}{Unif}

\DeclareMathOperator{\Var}{Var}

\DeclareMathOperator{\Tr}{Tr}

\DeclareMathOperator{\Id}{Id}
\DeclareMathOperator{\op}{op}

\newcommand{\calE}{\mathcal{E}}
\newcommand{\En}{\mathcal{E}_n}

\newcommand{\calU}{\mathcal{U}}

\newcommand{\calH}{\mathcal{H}}
\newcommand{\Sn}{\mathfrak{S}_n}

\newcommand{\tpi}{\widetilde{\pi}}

\newcommand{\norm}[1]{\left\lVert #1 \right\rVert}

\newcommand{\paren}[1]{\left( #1 \right)}

\newcommand{\set}[1]{\left\{ #1 \right\}}

\newcommand{\Indc}{\mathbf{1}}
\newcommand{\reals}{\mathbb{R}}

\newcommand{\rhobar}{\overline{\rho}}
\newcommand{\etabar}{\overline{\eta}}

\begin{document}

\title{Sharp Detection Threshold for Correlation among Multiple Unlabeled Gaussian Networks}

\author{
    Taha Ameen and Bruce Hajek 
    \thanks{
        The authors are with the Department of Electrical and Computer Engineering and the Coordinated Science Laboratory, University of Illinois, Urbana, IL 61801, USA. E-mails: \{ \texttt{tahaa3, b-hajek} \} \texttt{@illinois.edu} 
    }
}

\date{}

\maketitle

\begin{abstract}
    This paper studies the hypothesis testing problem of deciding whether $m\geq 2$ complete weighted graphs with Gaussian edge weights are mutually correlated after unknown relabelings of their vertices.  Under the null model all edge weights are independent standard Gaussians, whereas under the planted model the graphs share a latent vertex alignment and each pair of corresponding edge weights has correlation $\rho$.  For fixed $m$, we identify the sharp information-theoretic threshold for detection. Above the threshold, a generalized likelihood-ratio test achieves strong detection, whereas even weak detection is impossible below the threshold.  The result extends the two-graph detection threshold of Wu, Xu, and Yu~\cite{wu2023testing} to any fixed number of graphs, exhibits a side-information regime in which two graphs alone are insufficient but multiple graphs enable detection, and, together with the recovery threshold of Vassaux and Massouli\'e~\cite{vassaux2025}, shows that this Gaussian multi-graph model has no detection--recovery gap.
\end{abstract}

\tableofcontents

\section{Introduction}\label{sec:introduction}

Many modern data-analysis problems involve several observations of the same underlying collection of objects.  The observations may come from different platforms, laboratories, modalities, or time periods, and they are often correlated because they reflect a common ground truth.  At the same time, the labels attached to the individual objects need not agree across observations. For instance, users may be anonymized across social networks, proteins may be named differently across species, and vertices in two independently constructed graphs may have no common indexing convention.  Before one attempts to align such data, estimate a common structure, or transfer information from one source to another, a more basic question arises: is there statistically detectable evidence that the observations are correlated at all?

This question arises across several domains. For example, the network structure of two social networks such as Facebook and Twitter is correlated because users are likely to connect with the same individuals in both networks. This correlation can enable linkage or deanonymization attacks when the same users appear on different platforms~\cite{narayanan2008robust,narayanan2009deanonymizing}. A second example is from biology, where networks of interacting proteins may be correlated across species, indicating conserved functional structure across species~\cite{singh2008global}. Similarly, brain connectomes share anatomical organization across healthy subjects, and alignment can improve comparisons between individuals~\cite{sporns2005human,calissano2024graph}.  Related graph-matching formulations also appear in computer vision~\cite{schellewald2005probabilistic} and natural language processing~\cite{haghighi2005robust}.  These examples differ in their modeling details, but they share the same statistical obstruction: the correspondence between vertices is not observed.

This question of detecting correlation is naturally formulated as a hypothesis testing problem.  Under the null hypothesis the observed datasets are independent.  Under the alternative hypothesis they are correlated, but only after applying an unknown relabeling to each dataset.  The latent relabelings make the problem substantially harder than ordinary correlation testing: a test cannot directly compare the $i$-th coordinate across observations, because the coordinates have been scrambled.  Thus any successful procedure must find permutation-invariant evidence of correlation.

The present paper studies this problem for network-valued data.  We observe $m\geq 2$ weighted complete graphs on $n$ vertices.  Under the null model all edge weights in all graphs are independent standard Gaussians.  Under the planted model the graphs are aligned by an unknown tuple of vertex permutations, and the $m$ edge weights associated with each latent edge form an equicorrelated Gaussian vector with pairwise correlation $\rho$.  The central question is to determine the smallest order of $\rho$ for which the planted and null models can be distinguished.

The case of two graphs already captures the difficulty caused by unknown labels.  For $m=2$, Wu, Xu, and Yu~\cite{wu2023testing} identified the sharp Gaussian detection threshold, showing that correlation becomes detectable at scale $\rho^2 = 4 \log (n)/n$.  The multi-graph setting introduces a new phenomenon.  Each individual pair of graphs may contain too little correlation to be detected on its own, yet the aggregate evidence across many pairs can be detectable.  The question is therefore to establish how the detection threshold improves as additional correlated graphs are observed.

Our main result gives a sharp answer in the Gaussian model: for every fixed $m\geq2$, the detection threshold is
\[
        \rho_{\star}^2 = \frac{8}{m}\,\frac{\log n}{n}.
\]
Above this threshold a generalized likelihood-ratio test succeeds with vanishing total error; below this threshold the total variation distance between the planted and null distributions tends to zero, so no test can even weakly distinguish the two hypotheses.  The factor $1/m$ quantifies the gain from observing multiple graphs.  In particular, for each fixed $m>2$ there is a range of correlations in which two (in fact, any $m' < m$) graphs alone are information-theoretically insufficient, while the full collection of $m$ graphs contains enough evidence for reliable detection.

\subsection{Related literature}\label{sec:related-work}

\paragraph{Correlation detection for two graphs.}
The hypothesis testing problem for two unlabeled random graphs was initiated in a precise information-theoretic form by Wu, Xu, and Yu~\cite{wu2023testing}.  They established the sharp detection threshold for the Gaussian model, as well as the related Erd\H{o}s--R\'enyi setting, when the graphs are sufficiently dense.  Furthermore, in the sparse Erd\H{o}s--R\'enyi setting, they established the detection threshold up to a constant factor; the sharp constant was later obtained by Ding and Du~\cite{ding2023detection}.  In the constant-average-degree sparse regime, the corresponding strong-detection threshold was recently obtained by Feng~\cite{feng2025strong}.  A closely related multi-layer sparse Erd\H{o}s--R\'enyi detection problem was studied by Ochoa~\cite{ochoa2026detection}.  The tests attaining the information-theoretic thresholds are not known to run in polynomial time.  Consequently, computationally efficient tests and computational barriers have been studied, for instance by Barak, Chou, Lei, Schramm, and Sheng~\cite{barak2019nearly}, by Mao, Wu, Xu, and Yu~\cite{mao2024testing}, and by Ding, Du, and Li~\cite{ding2025low}.  Related computational transitions for correlated stochastic block models were studied by Chen, Ding, Gong, and Li~\cite{chen2026computational}.  Huang and Yang studied sampling variants of the Gaussian Wigner and Erd\H{o}s--R\'enyi models, as well as efficient tests based on bounded-degree motifs~\cite{huang2025samplecomplexity,huang2026graphsampling,huang2025motifs}.  More general weighted-graph distributions were considered by Oren-Loberman, Paslev, and Huleihel~\cite{oren2024testing}.  Related problems of correlation detection have also been investigated in other graph models, including Galton--Watson trees~\cite{ganassali2024correlation,ganassali2024statistical}, correlated graph growth models~\cite{racz2022correlated}, and correlated uniform attachment trees~\cite{baumler2026correlated}.

\paragraph{Graph alignment.}
Correlation detection is closely related to graph alignment, where the graphs are assumed to be correlated and the objective is to recover the hidden vertex correspondence.  For complete Gaussian graphs with $m=2$, the information-theoretic recovery threshold was independently established by Ganassali~\cite{ganassali2022sharp} and by Wu, Xu, and Yu~\cite{wu2022settling}.  Vassaux and Massouli\'e~\cite{vassaux2025} recently identified the recovery threshold for fixed $m\geq2$ in the multi-graph Gaussian model, and Even and Ganassali~\cite{even2025statistical} studied the corresponding Gaussian alignment problem when the number of graphs is allowed to grow with $n$, including low-degree evidence for computational barriers.  Algorithmic aspects of Gaussian graph alignment have also been studied through convex relaxations and approximate-message-passing methods~\cite{varma2025birkhoff,li2024robust}.  For correlated Erd\H{o}s--R\'enyi graphs, exact recovery was studied by Cullina and Kiyavash~\cite{cullina2016improved,cullina2017exact}; subsequent work addressed almost-exact and partial recovery~\cite{cullina2019kcore,wu2022settling,ding2023densesubgraph,hall2023partial,du2025optimal}, heterogeneous graphs~\cite{racz2023matching}, robustness to vertex corruptions~\cite{ameen2023robust}, partially correlated graphs~\cite{huang2025information}, and settings with more than two graphs~\cite{josephs2021recovery,ameen2024exact,ameen2024aligning}.  Correlated stochastic block models form another closely related direction: exact graph matching and exact community recovery were studied by R\'acz and Sridhar~\cite{racz2021correlated} and by Gaudio, R\'acz, and Sridhar~\cite{gaudio2022exact}, respectively, with subsequent work on efficient matching and on multiple correlated stochastic block models~\cite{chai2024efficient,racz2024harnessing}.  Finally, alignment models with node or contextual side information have recently been investigated in correlated contextual stochastic block models and Gaussian-feature models~\cite{yarandi2026contextual,huang2026attributed,zhang2021attributed}.

\paragraph{Database alignment.}
A parallel line of work studies unlabeled databases rather than graphs.  There, one observes collections of high-dimensional feature vectors that are independent under the null hypothesis and correlated through a latent permutation under the alternative.  Gaussian database models, including local decision-making tests, were studied by Tamir~\cite{tamir2022correlation,tamir2023correlation,tamir2025testing}; detection thresholds were obtained by Zeynep and Nazer~\cite{zeynep2022detecting}, sharpened by Elimelech and Huleihel~\cite{elimelech2023phase,elimelech2024detection} and by Jiao, Wu, and Xu~\cite{jiao_wu_xu2025}, and extended to general distributions by Paslev and Huleihel~\cite{paslev2024testing}.  Information-theoretic limits for recovering the latent database alignment have also been studied extensively~\cite{cullina2018fundamental,dai2019database,bakirtacs2021database,bakirtas2024database,dai2020achievability,dai2023gaussian}.  

\subsection{Contributions}\label{sec:contributions}

We consider $m\geq2$ complete Gaussian graphs on $n$ unlabeled vertices, with common pairwise correlation $\rho$ between corresponding edge weights.  The contributions of the paper are as follows.

First, we derive a sufficient condition for detection by thresholding the generalized likelihood ratio. Specifically, we show that strong detection is possible whenever
\[
        \rho^2 \geq \frac{8}{m}\,\frac{\log n}{n-1}.
\]
The proof combines a Gaussian quadratic-form concentration inequality with a union bound over all candidate alignments.  The leading term in the exponent is exactly balanced by the number $(n!)^{m-1}$ of possible alignments.
Second, we prove the matching converse: if
\[
        \rho^2 \leq \frac{8(1-\eta)}{m}\,\frac{\log n}{n}
\]
for a fixed $\eta\in(0,1)$, then the planted and null marginal distributions have total variation distance tending to zero.  The proof uses a truncated likelihood ratio in the spirit of~\cite{wu2023testing}.  The truncation conditions on the planted model to remove rare alignments with unusually large partial correlation sums.  After truncation, the second moment is controlled by a detailed trace expansion and by a sharp estimate from~\cite{vassaux2025}, for the number of alignments with prescribed overlaps.

The detection threshold obtained here agrees with the multi-graph Gaussian recovery threshold of~\cite{vassaux2025}.  Thus, in this model, there is no information-theoretic detection--recovery gap: the correlations at which one can reliably decide that a latent alignment exists are the same as the correlations at which one can recover that alignment.

\paragraph{Organization.}
\Cref{sec:problem} gives the formal model and definitions of strong and weak detection.  \Cref{sec:results} introduces the test statistic and states the main theorems.  \Cref{sec:achievability} proves the achievability result.  \Cref{sec:converse} proves the converse.  The derivation of the generalized likelihood-ratio statistic and other auxiliary estimates are collected in the appendices.

\section{Problem Formulation}\label{sec:problem}

Fix an integer $m\geq 2$ independent of $n$. For any finite set $A$, denote by $\binom{A}{2}\coloneqq\{B\subset A:|B|=2\}$ the collection of unordered pairs of distinct elements of $A$, and set
\[
        \En = \binom{[n]}{2},
        \qquad
        N = |\En| = \binom{n}{2} .
\]
We observe $m$ weighted complete graphs on the common vertex set $[n]$, encoded by their edge weights
\[
        X = \big(X^k_e\big)_{1\leq k\leq m,\, e\in\En}
        \in \R^{mN} .
\]
Let $\Sn$ denote the set of permutations of $[n].$  A permutation $\tau \in \Sn$ induces a permutation on $\En$, also denoted by $\tau$, defined by $\tau(\{i,j\}) = \{\tau(i),\tau(j)\}$ for $\{i,j\}\in \En.$  An {\em alignment} $\pi$ is an $m$-tuple $\pi=(\pi_1, \ldots , \pi_m) \in \Sn^m$ such that $\pi_1 = \Id.$
Under the \emph{null law} $\Q$, all $mN$ coordinates of $X$ are independent standard Gaussians. Under the \emph{planted law} $\P$, there is a latent alignment
\[
        \pi^*=(\pi^*_1,\ldots,\pi^*_m),
        \qquad
        \pi^*_1=\Id,
        \qquad
        \pi^*_2,\ldots,\pi^*_m \stackrel{\mathrm{iid}}{\sim}\Unif(\Sn),
\]
Conditionally on $\pi^*$, the edge tuples
\[
        \big(X^1_{\pi^*_1(e)},X^2_{\pi^*_2(e)},\ldots,X^m_{\pi^*_m(e)}\big) ,
        \qquad e\in\En,
\]
are independent centered Gaussian vectors with common covariance matrix
\begin{align} \label{eq:small_sigma}
        \Sigma=(1-\rho)\, I_m+\rho\, \mathbf{1}\mathbf{1}^{\top} ,
\end{align}
where $\rho\in[0,1)$ is the pairwise correlation between aligned edge weights. The data $X$ is observed; the latent alignment $\pi^*$ is not. The goal is to determine which hypothesis is true based on the observation. 

\paragraph{Additional notation}
We use $X^k$ as a shorthand for the collection of edge weights $(X^k_e)_{e\in\En}$ in the $k$-th graph. Given an alignment $\pi = (\pi_1, \ldots, \pi_m)$, write
\[
        \pi_{k\ell} \coloneqq \pi_\ell\circ \pi_k^{-1} ,
\]
for $k,\ell\in[m]$ so that $\pi_{k\ell}$ encodes the latent correspondence between the vertex labels of $X^k$ and $X^\ell$. Throughout, $\Q$ denotes the law of $X$ under the null hypothesis $H_0$, while $\P$ denotes the joint law of $(X,\pi^*)$ under the planted hypothesis $H_1$. With a slight abuse of notation, we write $\TV(\P,\Q)$ for the total variation distance between the marginal law of $X$ under $\P$ and the null law $\Q$. Standard asymptotic notation ($O(\cdot), o(\cdot), \Omega(\cdot), \cdots$) is used throughout, and it is implicit that $n\to\infty$. Further, we say $a_n \sim b_n$ if $a_n = (1+o(1)) \cdot b_n$.

\paragraph{Detection.}
A \emph{test} is a measurable function $\psi:\R^{mN}\to\{0,1\}$. The test $\psi$ achieves
\begin{itemize}
    \item \emph{strong detection} if its total error vanishes asymptotically:
    \[
        \P\bigl(\psi(X)=0\bigr) + \Q\bigl(\psi(X)=1\bigr) = o(1) ,
    \]
    \item \emph{weak detection} if it strictly outperforms random guessing:
    \[
        \P\bigl(\psi(X)=0\bigr) + \Q\bigl(\psi(X)=1\bigr) = 1 - \Omega(1) .
    \]
\end{itemize}
A standard fact from binary hypothesis testing is that strong detection is achievable if and only if $\TV(\P,\Q)=1-o(1)$, and weak detection is achievable if and only if $\TV(\P,\Q)=\Omega(1)$.

\medskip 



\section{Main Results}\label{sec:results}

For a candidate alignment $\pi$ set
\begin{equation}\label{eq:T-sigma}
        T_\pi(X)
        \coloneqq
        \sum_{e\in\En}\sum_{1\leq k<\ell\leq m}
        X^k_{\pi_k(e)}X^\ell_{\pi_\ell(e)} .
\end{equation}
Our analysis is based on the test statistic
\begin{equation}\label{eq:T-def}
        T(X)\coloneqq
        \max_{\pi_2,\ldots,\pi_m\in\Sn} T_\pi(X) ,
\end{equation}
which aggregates the aligned pairwise correlations over all graph pairs and over all edges, maximized over candidate alignments. The choice of $T$ is motivated by the fact that thresholding $T$ is approximately equivalent to the generalized likelihood-ratio test (GLRT) for distinguishing $\Q$ from $\P$; see \Cref{app:GLRT-derivation}.

The first result gives the achievability half of the detection threshold.

\medskip 

\begin{theorem}[Achievability]\label{thm:positive}
Suppose that
\begin{equation}\label{eq:positive}
        \rho^2 \geq \frac{8}{m}\cdot\frac{\log n}{n-1} \, .
\end{equation}
Then there exists a threshold $\tau$ such that the generalized likelihood-ratio test based on $T$ achieves strong detection, i.e.,
\(
    \P(T<\tau)+\Q(T\geq \tau)=o(1).
\)
\end{theorem}

The next result is a matching converse bound.

\medskip 

\begin{theorem}[Converse]\label{thm:negative}
Fix $\eta\in(0,1)$, and suppose that
\begin{equation}\label{eq:rho-bound-main}
        \rho^2 \leq \frac{8(1-\eta)}{m}\cdot\frac{\log n}{n} \, .
\end{equation}
Then $\TV(\P,\Q)\to 0$ as $n\to\infty$. In particular, weak detection is impossible in this regime.
\end{theorem}

Together, \Cref{thm:positive,thm:negative} identify the sharp critical scale
\[
    \rho_{\star}^2 \, = \, \frac{8}{m}\cdot\frac{\log n}{n} \, .
\]

The preliminary work~\cite{ameen2025detecting} established the achievability bound of \Cref{thm:positive} and a converse of the form
\[
    \rho^2 \, \leq \, \paren{\frac{4}{m-1} - \eps}\cdot \frac{\log n}{n} ,
\]
which differs from~\eqref{eq:rho-bound-main} by a multiplicative factor of $2(m-1)/m$. This bound was obtained by reducing the $m$-graph converse to the two-graph converse of~\cite{wu2023testing} via a genie-aided construction, which loses a constant factor at each inductive step. The present paper instead works directly with the likelihood ratio for $m$ graphs. The key technical ingredient is a second-moment computation applied to a suitably truncated likelihood ratio in the spirit of~\cite{wu2023testing}: the truncation handles the rare but large values of the likelihood ratio that would otherwise dominate the second moment. The argument also leverages a permutation-counting estimate due to Vassaux and Massouli\'e~\cite{vassaux2025} (recorded as \Cref{lemma:VM} below), which is  tight enough to close the multiplicative gap.

The paper~\cite{vassaux2025} establishes a sharp information-theoretic threshold for the closely related \emph{recovery} problem: given a sample from the planted law for $m$ Gaussian graphs, the latent alignment $\pi^*$ can be recovered with high probability if and only if $\rho^2 \geq (8/m)\cdot(\log n)/n$, up to lower-order terms. This recovery threshold and the detection threshold established here coincide. Consequently, the multi-graph Gaussian correlation problem exhibits no \emph{detection--recovery gap}: there is no region of correlations in which one can decide between the planted and null hypotheses but cannot recover the latent alignment. The two results are thus complementary, and together they give a complete information-theoretic picture of the $m$-graph Gaussian model.

\section{Achievability of detection}\label{sec:achievability}

The proof of \Cref{thm:positive} is a direct analysis of the statistic $T$ in~\eqref{eq:T-def}.  Under the planted law, the value of $T$ at the true alignment is a sum of independent Gaussian quadratic forms with positive mean.  Under the null law, a union bound over all alignments, combined with a one-sided Gaussian Hanson--Wright inequality, controls the maximum. Proofs  of the following two supporting lemmas are given in~\Cref{app:achievability-tools}.

Let $\|A\|_{\mathrm{op}}$ denote the spectral norm and $\|A\|_{\mathrm F}$ the Frobenius norm of a matrix $A.$

\medskip

\begin{lemma}[Gaussian Hanson--Wright bound~{\cite{hansonwright71}}]\label{lemma:HWGaussian}
Let $g\sim\mathcal N(0,I_d)$ and let $A$ be a real symmetric $d\times d$ matrix. For $Z=g^\top A g$ and every $t\ge0$,
\begin{align}
        \Prob\big(Z-\Expect Z\ge t\big)
        &\le
        \exp\left(
        -\frac{t^2}{4\left(\norm{A}_{\mathrm F}^2+\norm{A}_{\mathrm{op}}t\right)}
        \right).\label{eq:HWGaussian-one-sided}
\end{align}
Equivalently, for any fixed $\gamma\in(0,1)$,
\begin{align}
        \Prob\big(Z-\Expect Z\ge t\big)
        &\le
        \exp\left(
        -\frac14
        \min\left\{
        \frac{\gamma t^2}{\norm{A}_{\mathrm F}^2},
        \frac{(1-\gamma)t}{\norm{A}_{\mathrm{op}}}
        \right\}
        \right).\label{eq:HWGaussian-min-form}
\end{align}
\end{lemma}

The following lemma will be used to apply the Hanson-Wright bound in the context of this paper.
\medskip 

\begin{lemma}\label{lem:edge-quadratic-eigs}
Let
\(
        \Sigma_\rho=(1-\rho)I_m+\rho\mathbf 1\mathbf 1^\top
\)
and
\(
        Y\sim\mathcal N(0,\Sigma_\rho),
\)
and define
\(
        Z_\rho\coloneqq \sum_{1\le k<\ell\le m}Y_kY_\ell .
\)
Then $Z_\rho$ has the same distribution as $g^\top A_\rho g$, where $g\sim\mathcal N(0,I_m)$ and $A_\rho$ is a symmetric matrix with eigenvalues
\[
        \lambda_+(\rho)
        =\frac{m-1}{2}\big(1+(m-1)\rho\big)
        ~~~~\mbox{ and } ~~~~
        \lambda_-(\rho)
        =-\frac{1-\rho}{2}
\]
of multiplicity $1$ and $m-1$ respectively. 
\end{lemma}

The proof of~\Cref{thm:positive} consists of bounding the type-I and type-II error probabilities for the threshold defined by
\[
        \tau
        \coloneqq
        \binom{n}{2}\binom{m}{2}\rho - n^c,
\]
for a fixed constant $c\in(1,3/2)$.

\paragraph{Type-II error} Since $T\ge T_{\pi^*}$, it suffices to show that 
$\P(T_{\pi^*}<\tau) = o(1).$  Conditionally on $\pi^*$, the random variables
\[
        Z_e
        \coloneqq
        \sum_{1\le k<\ell\le m}
        X^k_{\pi^*_k(e)}X^\ell_{\pi^*_\ell(e)},
        \qquad e\in\En,
\]
are mutually independent, and each has the same distribution as
\[
        Z_\rho
        \coloneqq
        \sum_{1\le k<\ell\le m}Y_kY_\ell,
        \qquad
        Y\sim\mathcal N\left(0,(1-\rho)I_m+\rho\mathbf 1\mathbf 1^\top\right).
\]
By~\Cref{lem:edge-quadratic-eigs},
\(
        \Expect[Z_\rho]=\binom{m}{2}\rho
\)
and so we have
\(
    \Expect_\P[T_{\pi^*}]
    =\binom{n}{2}\binom{m}{2}\rho .
\)
Moreover, \Cref{lem:edge-quadratic-eigs} represents $T_{\pi^*}$ as a Gaussian quadratic form with
\[
        \norm{A}_{\mathrm{F}}^2
        =\binom{n}{2}
        \left(\lambda_+(\rho)^2+(m-1)\lambda_-(\rho)^2\right),
        \qquad
        \norm{A}_{\mathrm{op}}
        =\max\{|\lambda_+(\rho)|,|\lambda_-(\rho)|\},
\]
where
\[
        \lambda_+(\rho)
        =\frac{m-1}{2}\big(1+(m-1)\rho\big),
        \qquad
        \lambda_-(\rho)=-\frac{1-\rho}{2} .
\]
Since $m$ is fixed and $0\leq \rho<1$, these quantities are $O(n^2)$ and $O(1)$ respectively. 

Applying~\Cref{lemma:HWGaussian} to the negative of this quadratic form gives
\begin{align*}
        \P(T_{\pi^*}<\tau)
        &\leq
        \P\left(T_{\pi^*}-\Expect_\P[T_{\pi^*}]\leq -n^c\right) \\
        &\leq
        \exp\left\{
        -\frac{1}{4}
        \frac{n^{2c}}
        {\binom{n}{2}\left(\lambda_+(\rho)^2+(m-1)\lambda_-(\rho)^2\right)
        +\max\{|\lambda_+(\rho)|,|\lambda_-(\rho)|\}n^c}
        \right\} 
        =o(1) \, ,
\end{align*}
where the last step follows because $c>1$.

\paragraph{Type-I error} Under $\Q$, the distribution of $T$ does not depend on $\rho$, while the threshold $\tau$ is increasing in $\rho$. Hence the worst case for the bound is
\[
        \rho=\rho_n
        \coloneqq
        \sqrt{\frac{8\log n}{m(n-1)}} .
\]
For this value of $\rho$, let
\[
        \mu
        \coloneqq
        \binom{n}{2}\binom{m}{2}\rho_n,
        \qquad
        \tau=\mu-n^c .
\]
For all sufficiently large $n$, we have $\tau>0$. By the union bound over the $(n!)^{m-1}$ candidate alignments and by the permutation invariance of $\Q$,
\begin{align*}
        \Q(T\ge \tau)
        &\leq
        (n!)^{m-1} \, 
        \Q\left(
        \sum_{e\in\En}\sum_{1\le k<\ell\le m}X^k_eX^\ell_e
        \geq \tau
        \right).
\end{align*}
By \Cref{lem:edge-quadratic-eigs} with $\rho=0$, the random variable inside the last probability is a centered Gaussian quadratic form whose eigenvalues are $(m-1)/2$, with multiplicity $\binom{n}{2}$, and $-1/2$, with multiplicity $(m-1)\binom{n}{2}$. Therefore \Cref{lemma:HWGaussian} and the bound
\(
    n!\leq \exp\left((n+1/2)\log n-(n-1)\right)
\)
yield
\begin{align}\label{eq:ach-typeI-bound}
        \Q(T\ge\tau)
        \le
        \exp\left\{
        (m-1)(n+1/2)\log n-(m-1)(n-1)
        -\frac{\tau^2}
        {\binom{n}{2}m(m-1)+2(m-1)\tau}
        \right\}.
\end{align}
The choice of $\rho_n$ gives
\[
        \frac{\mu^2}{\binom{n}{2}m(m-1)}=(m-1)n\log n .
\]
Also, since $\mu=\Theta_m(n^{3/2}\sqrt{\log n})$ and $c<3/2$, we have $n^c=o(\mu)$. Using $1/(1+x)\ge 1-x$ for $x\ge0$ sufficiently small and $\tau^2\ge \mu^2-2\mu n^c$, we obtain
\[
        \frac{\tau^2}
        {\binom{n}{2}m(m-1)+2(m-1)\tau}
        \geq
        (m-1)n\log n
        -O\left(n^{c-1/2}\sqrt{\log n}\right)
        -O\left(n^{1/2}(\log n)^{3/2}\right).
\]
Both error terms are $o(n)$. Substituting this lower bound into \eqref{eq:ach-typeI-bound} gives
\[
        \Q(T\ge\tau)
        \le
        \exp\left
        \{
        -(m-1)(n-1)+\frac{m-1}{2}\log n+o(n)
        \right\}
        =o(1).
\]
Combining the two estimates proves
\(
        \P(T<\tau)+\Q(T\ge\tau)=o(1),
\)
as claimed.

\section{Impossibility of detection}\label{sec:converse}

The proof of \Cref{thm:negative} is by the second moment method. Throughout this section, $L=L(X)$ denotes the likelihood ratio between the marginal planted law $\P_X$ and the null law $\Q$. Since $\P_X\ll\Q$, Cauchy--Schwarz yields
\[
        \TV(\P,\Q)
        \, = \,
        \tfrac12\,\Expect_\Q |L-1|
        \, \leq \,
        \tfrac12\,\sqrt{\Var_\Q(L)} \, ,
\]
so that it suffices to show $\Expect_\Q[L^2]=1+o(1)$. As we discuss in \Cref{sec:converse-truncation}, however, the naive second-moment computation is dominated by atypical values of $L$ that, although rare under $\P$, contribute disproportionately to $\Expect_\Q[L^2]$. We therefore work with a truncated likelihood ratio that corresponds to an auxiliary measure $\P'$ close to $\P$ in total variation. The precise condition~\eqref{eq:rho-bound-main} on $\rho^2$ is only invoked near the end of the argument; weaker hypotheses on the rate at which $\rho\to 0$ are introduced as they become relevant.

The remainder of this section is organized as follows. \Cref{sec:converse-LR} sets up the matrix and permutation notation needed for the analysis and writes the likelihood ratio as a Gaussian integral indexed by alignments. \Cref{sec:converse-truncation} introduces the truncated planted law $\P'$ and computes the leading-order contribution to the second moment. \Cref{sec:converse-H} introduces a class $\calH$ of functions of permutations, that are exponentially negligible after averaging and absorbs the lower-order error terms in the second-moment expansion. \Cref{sec:converse-final} combines these ingredients to complete the proof.

\subsection{The likelihood ratio}\label{sec:converse-LR}

For a vertex permutation $\gamma\in\Sn$, write $\gamma(e)$ for its induced action on an unordered edge $e\in\En$. For an alignment $\pi=(\pi_1,\ldots,\pi_m)$, define the symmetric $mN\times mN$ matrix $K_\pi$ as follows. Its rows and columns are indexed by pairs $(k,e)\in[m]\times\En$, and
\begin{equation}\label{eq:K-entry}
        (K_\pi)_{(k,e),(\ell,f)}
        =
        \ind{k\ne \ell}\cdot
        \ind{\pi_{k\ell}(e)=f} .
\end{equation}
Equivalently, $(K_\pi)_{(k,e),(\ell,f)}=1$ if $k\ne\ell$ and the edge $f$ in graph $\ell$ corresponds to the same latent edge as $e$ in graph $k$.

Under $\P$ the covariance matrix of $X$ is  $\Sigma_{\pi^*}$ where
\[
        \Sigma_\pi = I_{mN}+\rho K_\pi,
\]
where $I_{mN}$ is the $mN\times mN$ identity matrix. We often write $I$ when the dimension is clear. Therefore the conditional likelihood ratio
\[
        L_\pi(X) \coloneqq \frac{\dd\P(X\mid\pi)}{\dd\Q(X)}
\]
is given by
\begin{equation}\label{eq:LLR}
        L_\pi(X)
        =
        \frac{1}{\sqrt{\det\Sigma_\pi}}
        \exp\left\{
        -\frac12 X^\top(\Sigma_\pi^{-1}-I)X
        \right\} .
\end{equation}
The full likelihood ratio is
\[
        L(X)=\Expect_\pi\big[L_\pi(X)\big].
\]
Hence
\[
        \Expect_\Q[L(X)^2]
        =
        \Expect_{\pi\ci\tpi}
        \left[
        \Expect_\Q\big[L_\pi(X)L_{\tpi}(X)\big]
        \right],
\]
where $\pi$ and $\tpi$ are independent alignments. 

We close this section with some notations and properties of permutations.
Given alignments $\pi$ and $\tpi$ let $d_{k\ell}(\pi, \tpi)$ denote the number of fixed points of $\pi_{\ell k}\circ\tpi_{k\ell}.$  It is equal to $|\{i\in [n]: \pi_{k\ell}(i) = \tpi_{k\ell}(i)\}|.$
Also let $\sigma$ be the alignment such that $\sigma_{1} = \mathrm{Id}$ and $\sigma_{k} = \pi_{k}^{-1}\circ\tpi_{k}$ for $2\leq k \leq m.$  
Note that 
\[
    \sigma_{k\ell} := \sigma_{\ell} \circ \sigma_k^{-1} = \pi_{\ell}^{-1}\circ \tpi_{\ell} \circ \tpi_k^{-1} \circ \pi_k = \pi_{\ell}^{-1}\circ \tpi_{k\ell} \circ \pi_k.
\]
Hence, $i$ is a fixed point of $\sigma_{k\ell}$ if and only if any of the following equivalent conditions hold:
\begin{align*}
    i =  \pi_{\ell}^{-1}\circ \tpi_{k\ell} \circ \pi_k (i) 
    ~~ \mbox{or} ~~
    \pi_{\ell}(i)=  \tpi_{k\ell} \circ \pi_k (i) 
    ~~ \mbox{or} ~~
    \pi_{\ell}\circ \pi_k^{-1}\circ \pi_k(i)=  \tpi_{k\ell} \circ \pi_k (i) \\
    ~~ \mbox{or} ~~
    \pi_{k\ell} \circ \pi_k(i) =  \tpi_{k\ell} \circ \pi_k (i) 
    ~~ \mbox{or} ~~
    \pi_k(i) = \pi_{\ell k}\circ \tpi_{k\ell} \circ \pi_k (i) .~~~~~~~~~~~~~~~~~~~~
\end{align*}
In other words, $i$ is a fixed point of $\sigma_{k\ell}$ if and only if $\pi_k(i)$
is a fixed point of $\pi_{\ell k}\circ\tpi_{k\ell}.$    Therefore, $d_{k\ell}(\pi,\tpi) = d_{k\ell}(\sigma),$ where  $ d_{k\ell}(\sigma)$ is defined to be the number of fixed points of $\sigma_{k\ell}.$    Therefore the following result from \cite{vassaux2025} is relevant and will play an important role near the end of the proof of \Cref{thm:negative}.

\medskip 

\begin{lemma} [\cite{vassaux2025}]  \label{lemma:VM}
     Let $(d_{k\ell})_{1 \leq k < \ell \leq m}$ be nonnegative integers.  Set
    \begin{align}
        F \big( (d_{k \ell}) \big) = \bigg| \Big\{\sigma \in \mathcal{S}_n^{m-1} \colon \,  d_{k\ell}(\sigma) \geq d_{k\ell} ~\mbox{ for all }~  1 \leq k < \ell \leq m, \Big\}\bigg| \, .
    \end{align}
    Then,\footnote{Here, the exponent is $ 2/m$ because the product is over $1\leq k < \ell \leq m$, instead of $1\leq k \neq \ell \leq m$ as in~\cite{vassaux2025}.  } 
    \begin{align} \label{usablef}
        F((d_{k\ell})) \leq \frac{(n!)^{m-1}}{\left(\prod_{1 \leq k  <  \ell \leq m} (d_{k\ell}!) \right)^{2/m}} \, .
    \end{align}
\end{lemma}
The same considerations as above hold for the edge permutations induced by the vertex permutations. (By definition, the edge permutation induced by a permutation $\tau\in\Sn,$ also called $\tau$, is defined by $\tau(e)=f$ if and only if $e=\{i,j\}$ for some $1\leq i < j \leq n$ and $f=\{\tau(i),\tau(j)\}.$)  
So similarly, let $D_{k\ell}(\pi,\tpi) = |\{ e\in \calE_n : \pi_{k\ell}(e) = \tpi_{k\ell}(e)\}|.$ Then
$D_{k\ell}(\pi,\tpi)= D_{k\ell}(\sigma)$, where $D_{k\ell}(\sigma)$ is the number of fixed edges of the edge permutation induced by $\sigma_{k\ell}.$   We write $D_{k\ell}(\pi,\tpi)$ or $D_{k\ell}(\sigma)$ interchangeably.

\subsection{Truncating the planted measure}\label{sec:converse-truncation}
We require an upper bound on $\Expect_{\Q}[ L_{\pi}(X)L_{\tpi}(X)]$ that is uniform enough to average over the independent alignments $\pi$ and $\tpi$.  The untruncated second moment is affected by rare values of $L_{\pi}(X)$ that are small-probability events under $\P$ but have a disproportionate effect under $\Q$.  We control these contributions by replacing $\P$ with a probability measure $\P'$ obtained by conditioning on a high-probability event.  Once we prove $\TV(\P',\Q)=o(1)$, the desired conclusion follows from the triangle inequality.
Specifically, $\P'$ is defined by setting
\[
    \P'(A) = \P(A|\calE) = \frac{ \P(A\cap \calE)}{\P(\calE)}
\]
for all Borel measurable $A \subset \reals^{m\binom n 2} \times\Sn^{m-1}$, where $\calE$ is defined as follows:
\[
    \calE
    \triangleq \bigcap_{1\le k < \ell \le m}\ \bigcap_{S\subset [n]:\, |S|\ge n/2} \calE_{k\ell,S} \, ,
\]
where for $1\leq k < \ell \leq m$ and $S\subset [n]$
\[
    \calE_{k \ell ,S}
    \triangleq 
    \bigg\{
        (X,\pi):
        \sum_{e\in \binom S 2}
        X^k_e \,X^{\ell}_{ \pi_{k\ell}(e) }
        \le
        \rho \binom{|S|}{2}+Cn^{3/2}
    \bigg\},
\]
where \(C>0\) is a sufficiently large constant depending only on \(m\).

It follows from~\cite[Lemma 1]{wu2023testing} that
\[
    \P \bigg( \, \bigcap_{S\subset [n]:\, |S|\ge n/2} \calE_{k\ell,S} \bigg) = 1 - e^{- \Omega(n)}.
\]
for $1\leq k < \ell \leq m.$   Since $\calE$ is the intersection of $\binom m 2$ such events, the following lemma holds:
\medskip
\begin{lemma}
    $\P(\calE) = 1 - e^{-\Omega(n)}.$
\end{lemma}

Recall for each $(k,\ell)$ with $1\leq k< \ell \leq m$ that
$d_{k\ell}(\pi,\tpi) = d_{k\ell}(\sigma) = 
| \{i\in [n] : \pi_{k\ell}(i)=\tpi_{k\ell}(i)\}|.$ 
For each ordered pair $k,\ell\in[m]$ with $k\ne\ell$, let
\[
    O_{k\ell} = 
    \begin{cases}
        \emptyset & \mbox{if } d_{k\ell}(\sigma)<  n / 2  \\        
        \{i\in [n] : \pi_{k\ell}(i)=\tpi_{k\ell}(i)\} & \mbox{if } d_{k\ell}(\sigma) \geq n / 2  
    \end{cases}
    \, .
\]
Note that 
\( 
    |O_{k \ell}|=d_{k\ell}(\sigma)\Indc_{\{d_{k\ell}(\sigma) \geq \frac n 2\}}. 
\)

We define a matrix $K'_{\pi}$ related to $K_{\pi}$ as follows:
\begin{align}\label{eq:K_prime-entry}
        (K'_\pi)_{(k,e),(\ell,f)}
        =
        \ind{k\ne \ell}\,
        \ind{\pi_{k\ell}(e)=f} \ind{e\in \binom{O_{k\ell}}{2}}
\end{align}
Note that $K_{\pi}$ depends only on $\pi$ while the sets $(O_{k\ell})$ and the matrix $K'_{\pi}$ depend on both $\pi$ and $\tpi$.
We collect here some matrix identities and trace expansions that will be used later. Their proofs are deferred to \Cref{app:matrix-properties}.
Recall that
\begin{align*}
        D_{k\ell}(\pi,\tpi) &
        \coloneqq
        \left|
        \set{e\in\En: \pi_{k\ell}(e)=\tpi_{k\ell}(e)} 
        \right| =D_{k\ell}(\sigma).
\end{align*}

\begin{lemma} \label{lemma:K_properties}
\begin{enumerate}[label=$\mathrm{(\alph*)}$]
    \item If $K'_{\tpi}$ is defined by \eqref{eq:K_prime-entry} with $\pi$ replaced by $\tpi,$ then $K'_{\pi}=K'_{\tpi}.$
    \item $K'_{\pi} \leq K_{\pi}$ pointwise and $K'_{\pi} \leq K_{\tpi}$ pointwise.
    \item $K_\pi$ and $K'_{\pi}$ are symmetric matrices.
    \item $  
    \begin{aligned} 
            \Tr(K_\pi K_{\tpi}) =  2\sum\nolimits_{1\le k<\ell\le m}D_{k\ell}(\sigma).
    \end{aligned}
    $
    \item
    \begin{align*} 
    \Tr(K_{\pi}K'_{\pi}) &= \Tr(K'_{\pi}K_{\pi})= \Tr(K'_{\pi}K_{\tpi}) = \Tr(K_{\tpi}K'_{\pi}) \\
    & = \Tr((K'_{\pi})^2) = 2\sum_{1\leq k < \ell \leq m}  \binom {d_{k\ell}(\sigma)} 2  \Indc_{\{ d_{k\ell}(\sigma)\geq \frac n 2\}}.
    \end{align*}
\end{enumerate}
\end{lemma}

We also record three higher order trace expansions that will be useful later in the converse proof.

\medskip 

\begin{lemma}[Trace expansions]\label{lem:trace-expansions}
For any two alignments $\pi$ and $\tpi$,
\begin{enumerate}[label=$\mathrm{(\alph*)}$]
    \item \label{trace-expansion-W3}
    $\displaystyle
    \begin{aligned}[t]
    W_3(\pi,\tpi) \coloneqq \Tr(K_\pi^2K_{\tpi})
    &=
    \sum_{\substack{k,\ell,k'\in[m]\\
            k,\ell,k'\,\mathrm{distinct}}}
    \Big|
    \left\{
            e\in\En:
            \big(
            \tpi_{k'k}
            \circ \pi_{\ell k'}
            \circ \pi_{k\ell}
            \big)(e)=e
    \right\}
    \Big| .
    \end{aligned}
    $
    \item \label{trace-expansion-W4}
    $\displaystyle
    \begin{aligned}[t]
    W_4(\pi,\tpi) \coloneqq \Tr\big((K_\pi K_{\tpi})^2\big)
    &=
    \sum_{\substack{k,\ell,k',\ell'\in[m]\\
                     \{k,k'\}\cap \{\ell, \ell'\} = \emptyset}}
    \left|
    \left\{
            e\in\En:
            \big(
            \tpi_{\ell'k}
            \circ \pi_{k'\ell'}
            \circ \tpi_{\ell k'}
            \circ \pi_{k\ell}
            \big)(e)=e
    \right\}
    \right| .
    \end{aligned}
    $
    \item \label{trace-expansion-tildeW4}
    $\displaystyle
    \begin{aligned}[t]
    \widetilde{W}_4(\pi,\tpi) \coloneqq \Tr(K_\pi^3K_{\tpi})
    &=
    \sum_{\substack{k,\ell,k',\ell'\in[m]\\
            \{k,k'\} \cap \{\ell,\ell'\} = \emptyset }}
    \Big|
    \left\{
            e\in\En:
            \big(
            \tpi_{\ell'k}
            \circ \pi_{k'\ell'}
            \circ \pi_{\ell k'}
            \circ \pi_{k\ell}
            \big)(e)=e
    \right\}
    \Big| .
    \end{aligned}
    $
\end{enumerate}
\end{lemma}

For $1\leq k < \ell \leq m,$ if $d_{k\ell}(\sigma) < n/2$ then $O_{k\ell} = \emptyset$, and if $d_{k\ell}(\sigma) \geq n/2$ then $\calE_{k\ell ,O_{k\ell}} \subset \cal E.$ 
Therefore, on the event $(X,\pi)\in \calE$, 
\begin{align*}
    \rho X^{\top}K'_{\pi}X 
    \, = \, 
    2 \rho \sum_{1 \leq k < \ell \leq m} \sum_{e \in \binom{O_{k\ell}}{2}} X^k_e X^{\ell}_{\pi_{k\ell}(e)}
    \, \leq \,  \tau_{\sigma}
\end{align*}
where
\begin{align*}
   \tau_{\sigma} &= 2 \sum_{1\leq k < \ell \leq m} \left\{ \rho^2 \binom {d_{k\ell}(\sigma)} 2 + \rho \, Cn^{3/2} \right\} \Indc_{\{ d_{k\ell}(\sigma)\geq \frac n 2\}} 
\end{align*}
Therefore, for every $r \geq 0$,
\[
    \Indc_{(X,\pi) \in \calE}
    \leq
    \exp\{r(\tau_\sigma-\rho X^{\top}K'_{\pi}X )\} \, .
\]
We will set $r=1/2$ below; for the moment assume only that $0\leq r\leq 1$.
Consider the conditioned likelihood ratio,
\[
    L'(X)
    =
    \frac{1}{\P(\calE)} \, 
    \Expect_{\pi}
    \left[
        L_{\pi}(X) \cdot \Indc_{(X,\pi) \in \calE}
    \right].
\]
Its second moment may be written as 
\begin{align}
    \Expect_{\Q}[(L')^2]
    &=
    \frac{1}{\P(\calE)^2} \, 
    \Expect_{\pi,\tpi}
    \Expect_{\Q}
    \left[
        L_{\pi}(X)L_{\tpi}(X) \cdot
        \Indc_{(X,\pi) \in \calE}
        \Indc_{(X,\tpi) \in \calE}
    \right] \nonumber
    \\
    & \leq
        \frac{1}{\P(\calE)^2} \, 
    \Expect_{\pi,\tpi}
    \Expect_{\Q}
    \left[
        L_{\pi}(X)L_{\tpi}(X) \cdot
        \Indc_{(X,\pi) \in \calE}
    \right]  \nonumber
    \\
    &\leq
    \frac{1}{\P(\calE)^2} \, 
    \Expect_{\pi,\tpi}
    \Expect_{\Q}
    \left[
        L_{\pi}(X)L_{\tpi}(X)
        \exp\{r(\tau_\sigma-\rho X^{\top}K'_{\pi}X )\}
    \right].  \label{eq:partial_progress}
\end{align}
Since $\P(\calE)=1-o(1)$, the prefactor $\P(\calE)^{-2}$ is $1+o(1)$.
Hence, a sufficient condition for $\TV(\P',\Q)=o(1)$ is 
\(
    \Expect_{\pi,\tpi}
    \Expect_{\Q}
    \left[
        L_{\pi}(X)L_{\tpi}(X)
        \exp\{r(\tau_\sigma-\rho X^{\top}K'_{\pi}X )\}
    \right] = 1 + o(1) \, .
\)
The following proposition will be used to bound the inner expectation in \eqref{eq:partial_progress}.  

\medskip

\begin{proposition}   \label{prop:det-id}  
    Suppose that $\rho^2 = o(1).$
    Given two alignments $\pi$ and $\tpi$ and $r\in [0,1]$, for sufficiently large $n$, 
    \begin{align}\label{eq:det-series}
        \log \Expect_\Q\Big[L_\pi(X)L_{\tpi}(X) \exp\{r(\tau_\sigma-\rho X^{\top}K'_{\pi}X )\}\Big]
        = 
        r\tau_{\sigma} + \frac12\sum_{s\ge 1}\frac 1 {s}
        \Tr\left[A^s\right] .
    \end{align}
    where
    \[
        A = \rho^2K_\pi K_{\tpi}  - 2r\rho K'_{\pi} - 2r\rho^2 (K_{\pi}K'_{\pi} + K'_{\pi}K_{\tpi} ) - 2r\rho^3   K_{\pi}K'_{\pi}K_{\tpi} .
    \]
\end{proposition}

\begin{proof}
Using \eqref{eq:LLR},
\begin{align*}
    \Expect_\Q & \big[L_\pi(X)L_{\tpi}(X) \exp\{r(\tau_\sigma-\rho X^{\top} K'_{\pi}X )\}\big]  
    \\[0.5em]
    & =
    \frac{1}{\sqrt{\det\Sigma_\pi}\sqrt{\det\Sigma_{\tpi}}} \, 
    \Expect_\Q\left[
    \exp\left\{r\tau_{\sigma}
    -\frac12X^\top(\Sigma_\pi^{-1}+\Sigma_{\tpi}^{-1} - 2I  + 2r\rho K'_{\pi})X
    \right\}
    \right] \notag
    \\[0.5em]
    & \stackrel{\mathrm{(a)}}{=}
    \frac{e^{r\tau_{\sigma}}}{\sqrt{\big. \det\Sigma_\pi}\sqrt{\big. \det\Sigma_{\tpi}}\sqrt{\det(\Sigma_\pi^{-1}+\Sigma_{\tpi}^{-1} - I + 2r\rho K'_{\pi})}} \notag
    \\[0.5em]
    & \stackrel{\mathrm{(b)}}{=}
    e^{r\tau_{\sigma}}\det\big(\Sigma_\pi+\Sigma_{\tpi}-\Sigma_\pi\Sigma_{\tpi} + 2r\rho \Sigma_{\pi} K'_{\pi}\Sigma_{\tpi}\big)^{-1/2} 
    \\[0.5em]
    & \stackrel{\mathrm{(c)}}{=}
    e^{r\tau_{\sigma}}\det(I-A)^{-1/2},
\end{align*}
where step (a) uses the Gaussian quadratic-form identity
\[
        \Expect\big[ \, e^{-\frac12 x^\top Mx} \, \big]
        =
        \det(I+M)^{-1/2},
\]
valid for $x\sim N(0,I)$ whenever $I+M$ is positive definite. The identity is applicable here because the matrix in the determinant is $I+O(\rho)$ in operator norm, uniformly in $\pi,\tpi$, and is therefore positive definite for all sufficiently large $n$.  Step (b) follows by multiplying the matrix inside the determinant on the left by $\Sigma_\pi$ and on the right by $\Sigma_{\tpi}$, and using multiplicativity of the determinant.  Step (c) follows from $\Sigma_\pi=I+\rho K_\pi$ and $\Sigma_{\tpi}=I+\rho K_{\tpi}$.
Since each row of $K_\pi, K_{\tpi}$ and $K'_{\pi}$ has at most $m-1$ nonzero entries and they are all one, Gershgorin's circle theorem implies
\(
\max\{         \norm{K_\pi}_{\op},  \norm{K_{\tpi}}_{\op}, \norm{K'_{\pi}}_{\op}\} \le m-1 .
\)
Therefore
\begin{align} \label{eq:norm-A}
        \norm{A}_{\op} = O(\rho) = o(1)
\end{align}
so for all sufficiently large $n$, $\norm{A}_{\op}<1$, so
\(
        -\log(I-A)=\sum_{s\ge 1} {A^s}/{s} \, .
\)
Taking traces and using $\log\det B=\Tr(\log B)$ yields
\[
        -\log\det(I-A)
        =
        \sum_{s\ge 1}\frac{\Tr(A^s)}{s}. \qedhere
\]
\end{proof}

We next examine some of the terms appearing on the right-hand side of \eqref{eq:det-series}.
Using the fact $K'_{\pi}$ has zero diagonal and zero trace,~\Cref{lemma:K_properties}(e),
the fact $\Tr(AB)=\Tr(BA)$ for matrices $A,B$, and dropping negative terms yields
\begin{align*}   
    \frac 1 2 \Tr(A) & \, \leq \, \frac 1 2 \rho^2 \Tr\left(K_\pi K_{\tpi} - 4r (K'_{\pi})^2  \right),
    \\
    \frac 1 4 \Tr(A^2) & \, \leq  \, r^2\rho^2 \Tr((K'_\pi)^2)  + J_2 + o(\rho^5) \\
    \Tr(A^3) &\leq 3(2r)^2\rho^4 \Tr((K'_{\pi})^2K_{\pi}K_{\tpi}) + o(\rho^5)  \\
    \Tr(A^4) &\leq 16r^2\rho^4 \Tr((K'_{\pi})^4)  + o(\rho^5)  \, , 
\end{align*}
where
\begin{align*}
    J_2 &= 2 r^2\rho^3 \Tr \big(K_{\pi}K'_{\pi}K'_{\pi} +K'_{\pi} K'_{\pi}K_{\tpi}) 
    + \frac{\rho^4} 4 \Tr((K_{\pi} K_{\tpi})^2) \\
    &+ r^2\rho^4\Tr((K_{\pi}K'_{\pi} + K'_{\pi}K_{\tpi} )^2) + 2r^2\rho^4\Tr(K'_{\pi} K_{\pi}K'_{\pi}K_{\tpi})  \, .
\end{align*}

\begin{lemma} \label{lemma:bounding}
Suppose that $n^2\rho^5 = o(1).$
With $W_3,W_4,\widetilde{W}_4$ as defined in~\Cref{lem:trace-expansions}, the following bounds hold:
\begin{align*}
    J_2 & \, \leq \, 
        4\rho^3 W_3(\pi,\tpi) +  \frac {17\rho^4} 4 W_4(\pi,\tpi) + 2 \rho^4 \widetilde{W}_4(\pi,\tpi)  
    \\[0.4em]
    \Tr(A^3) & \, \leq \, 
        12 \rho^4 W_4(\pi,\tpi)  + o(1)  
    \\[0.5em]
    \Tr(A^4) & \, \leq \, 
        16 \rho^4 W_4(\pi,\tpi)  + o(1)
\end{align*}
Also, the sum of all terms $s\ge 5$ in \eqref{eq:det-series} is $O(N\rho^5)=o(1)$.
\end{lemma}
\begin{proof}
We use $r\leq 1.$
Recall that $K'_{\pi} \leq K_{\pi}$ and $K'_{\pi} \leq K_{\tpi}$ pointwise.  Thus, 
\[
    \Tr(K_\pi K'_\pi K'_\pi) \leq \Tr(K_\pi^2K_{\tpi})=W_3(\pi,\tpi),
    \qquad
    \Tr(K'_\pi K'_\pi K_{\tpi}) \leq \Tr(K_\pi^2K_{\tpi})=W_3(\pi,\tpi).
\]
Similarly, $K_{\pi}K'_{\pi} + K'_{\pi}K_{\tpi}  \leq 2K_{\pi}K_{\tpi}$ so
$\Tr((K_{\pi}K'_{\pi} + K'_{\pi}K_{\tpi} )^2)\leq 4\Tr((K_{\pi}K_{\tpi})^2)= 4 W_4(\pi,\tpi). $   Also, $K'_{\pi}K_{\pi}K'_{\pi}K_{\tpi} \leq K_{\pi}^3K_{\tpi}$
so $\Tr(K'_{\pi}K_{\pi}K'_{\pi}K_{\tpi}) \leq \widetilde{W}_4(\pi,\tpi).$
Finally, $\Tr((K'_{\pi})^4)\leq \Tr((K'_{\pi})^2K_{\pi}K_{\tpi})\leq W_4(\pi,\tpi)$.
These yield the bounds in the centered equations in the lemma.

Recall that $\|A\|_{\op} = O(\rho)$ uniformly in $n$. 
Hence $\Tr(A^s) \leq Nm\|A\|_{\op}^s \leq Nm(C\rho)^s$ for some constant $C,$ implying the last statement of the lemma.
\end{proof}

Assume for the remainder of this section that $n^2\rho^5\to 0$ and $n \rho^2 \to \infty$.
Applying ~\Cref{prop:det-id},  ~\Cref{lemma:K_properties} and~\Cref{lemma:bounding} and recalling
\begin{align*}
    r\tau_{\sigma} & \, = \,  2r \sum_{1\leq k < \ell \leq m}  \left\{ \rho^2 \binom {d_{k\ell}(\sigma)} 2 + \rho \, Cn^{3/2} \right\} \cdot \Indc_{\{ d_{k\ell}(\sigma)\geq \frac n 2\}}  
\end{align*} 
we have for a constant $C$
\begin{align*}
    \log \Expect_\Q\big[L_\pi(X)L_{\tpi}(X) \Indc_{(X,\pi) \in \calE}\big]
    &       \leq \rho^2(2r-4r+2r^2)  \sum_{1\leq k < \ell \leq m} \! \binom {d_{k\ell}(\sigma)} 2 \cdot \Indc_{\{d_{k\ell}(\sigma)\geq \frac n2\}}     +\rho^2 \!\!\!  \sum_{1\leq k < \ell \leq m} \!\!\! D_{k\ell}(\sigma) \\[0.5em]
    &~~~~~        + 2\rho r \, Cn^{3/2}  \sum_{1\leq k < \ell \leq m}  \Indc_{\{ d_{k\ell}(\sigma)\geq \frac n 2\}}  \\[0.5em]
    &~~~~  + C\rho^3 \,  W_3(\pi,\tpi)   + C\rho^4 \Big(W_4(\pi,\tpi)  + \widetilde{W}_4(\pi,\tpi) \Big)  + o(1) \, .
\end{align*}

Setting $r=1/2$ then yields
\begin{align}
    \log \Expect_\Q & \big[L_\pi(X)L_{\tpi}(X) \Indc_{(X,\pi) \in \calE}\big]
    \nonumber\\
    & 
    \leq   - \frac 1 2 \, \rho^2 \!\! \sum_{1\leq k < \ell \leq m} \binom {d_{k\ell}(\sigma)} 2 \Indc_{\{d_{k\ell}(\sigma)\geq \frac n2\}}
    \, + \,  \rho^2 \!\! \sum_{1\leq k < \ell \leq m}  D_{k\ell}(\sigma) \nonumber 
    \\
    &~~~~~        + \rho \, Cn^{3/2}  \sum_{1\leq k < \ell \leq m} \Indc_{\{ d_{k\ell}(\sigma)\geq \frac n 2\}}  
    + C\rho^3 \, W_3(\pi,\tpi) + C\rho^4 \Big(W_4(\pi,\tpi)  + \widetilde{W}_4(\pi,\tpi) \Big) + o(1)   \nonumber 
    \\[0.5em]
    &  \stackrel{\mathrm{(a)}}{=}  \frac 1 2 \sum_{1\leq k < \ell \leq m}  \rho^2 \binom {d_{k\ell}(\sigma)} 2  \Indc_{\{ d_{k\ell}(\sigma)\geq \frac n 2\}} 
    +\rho^2  \sum_{1\leq k < \ell \leq m}  \binom {d_{k\ell}(\sigma)} 2 \Indc_{\{d_{k\ell}(\sigma)<\frac n 2 \} }  \nonumber 
    \\
    &~~~~~     
    +  \rho \,  C n^{3/2} \sum_{1 \leq k < \ell \leq m} \Indc_{ \{ d_{k\ell}(\sigma) \geq \frac n2 \}} +  Y(\pi,\tpi)  \nonumber \\
    & \stackrel{\mathrm{(b)}}{\leq}  (1+o(1)) \, \rho^2 \sum_{1\leq k < \ell \leq m} 
    \frac{ 1 +  \Indc_{\{ d_{k\ell}(\sigma) <  \frac n 2\}} }  2 
    \binom {d_{k\ell}(\sigma)} 2 + Y(\pi,\tpi)   \label{eq:term_bnd_with_Y}
\end{align}
where
\begin{align*}
    Y(\pi,\tpi) 
    = \rho^2  \sum_{1\leq k < \ell \leq m}  
    \left[ D_{k\ell}(\sigma)- \binom {d_{k\ell}(\sigma)} 2 \right]
    + C\rho^3 W_3(\pi,\tpi) + C\rho^4 \Big( W_4(\pi,\tpi)  + \widetilde{W}_4(\pi,\tpi) \Big) + o(1) \, .
\end{align*}
Here, (a) decomposes
\[
    D_{k\ell}(\sigma)
    =
    \binom{d_{k\ell}(\sigma)}2
    +
    \left[
        D_{k\ell}(\sigma)-\binom{d_{k\ell}(\sigma)}2
    \right]
\]
and (b) uses that, on \(\{d_{k\ell}(\sigma)\ge n/2\}\), $\rho \, Cn^{3/2} = o\left(\rho^2\binom {d_{k\ell}(\sigma)} 2\right)$ by the assumption $n\rho^2 \to \infty$ (i.e. $\sqrt{n} \rho \to \infty).$

The next section shows that the term $Y(\pi,\tpi)$ can be neglected after the averaging over $(\pi,\tpi).$

\subsection{Exponentially negligible functions}\label{sec:converse-H}

The main result of this section is~\Cref{lemma:terms_in_H} below.
\medskip
\begin{definition}
    Let $\calH$ denote the set of nonnegative functions $Y_n(\pi,\tpi)$ on $(\Sn^{m-1})^2$  such that for every constant $\Gamma >0$, 
    \[
        \Expect_{\pi\perp\tpi} \big[\exp (\Gamma Y) \big] = 1 + o(1) ~~~~  \mbox{ as } n\to\infty \, .
    \]
\end{definition}

\medskip
\begin{lemma} \label{lemma:H_closure}
    \begin{enumerate}[label=$\mathrm{(\alph*)}$]
        \item If $Y_1, Y_2, \ldots , Y_J \in \calH$ for a constant number of terms $J$, then $Y_1 + \ldots +  Y_J \in \calH.$  
        \item If $Y=f(n)$ for some function $f$ such that $f=o(1)$, then $Y\in \calH.$
    \end{enumerate}
\end{lemma}
\begin{proof}  
By the inequality $x_1x_2\cdots x_J \leq (x_1^J + \cdots +  x_J^J)/ J$, 
\[
    \exp \big( \Gamma(Y_1 + \ldots +  Y_J ) \big) 
    \, \leq \, 
    \frac 1 J \big( \exp (J \, \Gamma Y_1) + \ldots +  \exp (J \, \Gamma Y_J) \big),
\] 
which implies part (a).  Part (b) follows from $\exp(\Gamma \cdot o(1)) = 1+o(1)$.  
\end{proof} 

The following lemma follows directly from classical results about permutations -- the proof is in the appendix.
\medskip

\begin{lemma} \label{lem:exp-Cj-square}
    Fix an integer $j\ge 1$.  Let $C_j=C_j^{(n)}$ denote the number of $j$--cycles in a uniformly random permutation of $\{1,\dots,n\}$.  
    \begin{enumerate}[label=$\mathrm{(\alph*)}$]
        \item If $\theta_n = o(1)$ then $ \mathbb{E}\big[e^{\theta_n C_j}\big] = 1 + o(1)$.
        \item If $\theta_n = o\big( \frac{\log n}{n}\big)$ then $\mathbb{E}\big[e^{\theta_n C_j^2}\big]=1 + o(1).$
    \end{enumerate}
\end{lemma}

\medskip

\begin{lemma}  \label{lemma:terms_in_H}
    Suppose $\frac{n \rho^4}{\log n}  \to 0$ and $\rho^3=o\big(\frac{\log n}{n}\big)$.
    Then  $ \rho^2  \sum_{1\leq k < \ell \leq m}  \left[ D_{k\ell}- \binom {d_{k\ell}(\sigma)} 2 \right]$,
    $\rho^3W_3(\pi,\tpi)$, $\rho^4W_4(\pi,\tpi)$, and $\rho^4 \widetilde{W}_4(\pi,\tpi)$ are in $\calH.$ 
\end{lemma}

\begin{proof}
By~\Cref{lemma:H_closure}(a), to show
$ \rho^2  \sum_{1\leq k < \ell \leq m}  \big[ D_{k\ell}- \binom {d_{k\ell}(\sigma)} 2 \big] \in \calH,$ it suffices
to show that the individual terms $\rho^2\big(D_{k\ell}- \binom {d_{k\ell}(\sigma)} 2\big)\in \calH$ for all $k, \ell$ with $1 \leq k < \ell \leq m.$

Now, $D_{k\ell}(\sigma) - \binom {d_{k\ell}(\sigma)} 2  = c_{k\ell}$ where
\[
    c_{k\ell}
    \coloneqq
    \big|
        \set{\{i,j\}\subset[n]:\sigma_{k\ell}(i)=j,\ \sigma_{k\ell}(j)=i}
    \big| \, .
\]
In other words $c_{k\ell}$ is the number of 2-cycles (transpositions) in the cycle decomposition of the vertex permutation $\sigma_{k\ell}.$
Since $\sigma_{k\ell}$ is uniformly distributed over $\Sn$, $\rho^2 c_{k\ell} \in \calH$ by~\Cref{lem:exp-Cj-square}.

\paragraph{Bounding $\rho^3W_3(\pi,\tpi)$.}
By ~\Cref{lem:trace-expansions}\textup{(a)}, it suffices by~\Cref{lemma:H_closure}(a) to show that each summand is in $\calH$ after multiplication by $\rho^3$.
Given distinct $k,\ell,k'\in[m]$, let
\[
    \tau = \tpi_{k'k}\circ \pi_{\ell k'}\circ \pi_{k\ell}.
\]
Since $\tpi_{k'k}$ is independent of the other two permutations in $\tau$ and is uniformly distributed over $\Sn$, $\tau$ is also uniformly distributed over $\Sn$. Therefore,
\begin{align}  \label{eq:fixE_tau_in_H}
    \rho^3\operatorname{FixE}(\tau)
    =
    \rho^3\binom{C_1(\tau)}2+\rho^3C_2(\tau) \in \calH,
\end{align}
where the first term uses $\rho^3=o(\log n/n)$ and~\Cref{lem:exp-Cj-square}\textup{(b)}, while the second uses \Cref{lem:exp-Cj-square}\textup{(a)}.
Hence $\rho^3W_3(\pi,\tpi)\in\calH$.

\paragraph{Bounding $\rho^4 W_4(\pi,\tpi)$.}
By \Cref{lem:trace-expansions}\textup{(b)}, 
$ W_4(\pi,\tpi)$ is the sum of
\[
    D_{k\ell k'\ell'} =
    \big|
    \left\{
            e\in\En:
            \big(
            \tpi_{\ell'k}
            \circ \pi_{k'\ell'}
            \circ \tpi_{\ell k'}
            \circ \pi_{k\ell}
            \big)(e)=e
    \right\}
    \big| 
\]
over $k,\ell,k',\ell' \in [m]$ such that
     $ \{k,k'\} \, \cap \, \{\ell,\ell' \} = \emptyset.$
By~\Cref{lemma:H_closure}(a) it is sufficient to show that
$\rho^4 D_{k\ell k'\ell'}
\in \calH$ for all such $k,\ell,k',\ell'.$ We consider two cases. 

First consider the case $k=k'$ and $\ell = \ell'$ corresponding to terms of the
form $D_{k\ell k\ell}.$   Note that
$D_{k\ell k\ell}$ is the number of edges fixed by $\tau^2$ where $\tau = \tpi_{\ell k}\circ\pi_{k\ell}$  and $\tau$ is uniformly distributed over $\Sn.$   A vertex is fixed by $\tau^2$ if it lies in a $1$-cycle or a $2$-cycle of $\tau$. Therefore the number of vertices fixed by $\tau^2$ is
\(
        C_1(\tau)+2 \, C_2(\tau).
\)
Edges between such vertices are fixed by $\tau^2$. In addition, each $4$-cycle of $\tau$ contributes two transpositions to $\tau^2$, and hence contributes two more fixed edges. Thus
\begin{align}
  \rho^4 D_{k\ell k\ell}  =\rho^4 \operatorname{FixE}(\tau^2)
   & =\rho^4
    \binom{C_1(\tau)+2 \, C_2(\tau)}{2}
    +
    2\rho^4 \, C_4(\tau) \nonumber \\
   & \leq 2\rho^4C^2_1(\tau) + 8\rho^4C^2_2(\tau) + 2\rho^4C_4(\tau)  \label{eq:D4_bnd}
\end{align}
By~\Cref{lem:exp-Cj-square} each term on the right-hand side of \eqref{eq:D4_bnd} is in $\calH$
so by~\Cref{lemma:H_closure}, $ \rho^4 D_{k\ell k\ell} \in \calH.$

Next, consider the second case, namely, that either $k\neq k'$ or $\ell \neq \ell'$ or both. We assume $k\neq k'$; the proof for $\ell\neq\ell'$ is similar.  Either $k\neq 1$ or $k'\neq 1.$   We assume $k\neq 1$; the case $k'\neq 1$ is similar.
Let $\tau = \tpi_{\ell'k} \circ \pi_{k'\ell'} \circ \tpi_{\ell k'} \circ \pi_{k\ell}.$
Using the definition $\pi_{k\ell}=\pi_{\ell}\pi_{k}^{-1}$ and similarly for $\tpi_{k\ell}$, we find
\[
\tau = \tpi_k\circ\tpi_{\ell'}^{-1}\circ\pi_{\ell'}\circ\pi_{k'}^{-1}\circ\tpi_{k'}\circ\tpi_{\ell}^{-1}\circ\pi_{\ell}\circ\pi_k^{-1}.
\]
Note that $\tpi_k$ appears only once and is uniformly distributed over $\Sn$ and is independent of the other permutations defining $\tau.$  Hence $\tau$ is also uniformly distributed over $\Sn.$  Thus, recalling \eqref{eq:fixE_tau_in_H}, $\rho^4 D_{k,\ell,k,\ell'} = \rho^4  \operatorname{FixE}(\tau) \leq \rho^3
\operatorname{FixE}(\tau) \in \calH.$  

The two cases cover all possibilities for $k,\ell,k',\ell'$, so $\rho^4 W_4(\pi,\tpi) \in \calH.$  

\paragraph{Bounding $\rho^4 \widetilde{W}_4(\pi,\tpi)$.} Note that $\widetilde{W}_4(\pi,\tpi)$ is the same as $W_4(\pi,\tpi)$ but with $D_{k\ell k' \ell'}$ being the number of fixed points of
$  \tpi_{\ell'k}  \circ \pi_{k'\ell'}  \circ \pi_{\ell k'}   \circ \pi_{k\ell}.$
Since $\tpi$ appears only once in each composition, all terms fall into the second case of the proof that $W_4(\pi,\tpi) \in \calH$ so that
$\rho^4\widetilde{W}_4(\pi,\tpi) \in \calH.$
\end{proof}

\subsection{Finishing the proof}\label{sec:converse-final}

Recall from \eqref{eq:term_bnd_with_Y} that if $n^2\rho^5\to 0$ and $n\rho^2\to \infty$ then
\begin{align*}
    \log \Expect_\Q\big[L_\pi(X)L_{\tpi}(X) \Indc_{(X,\pi) \in \calE}\big]
    & 
    \leq  \big(1+o(1) \big) \, \rho^2 \!\! \sum_{1\leq k < \ell \leq m} \frac{ 1 + \Indc_{\{ d_{k\ell}(\sigma) < \frac n 2\}}} 2 \binom {d_{k\ell}(\sigma)} 2  + Y(\pi,\tpi),
\end{align*}
where by~\Cref{lemma:H_closure} and~\Cref{lemma:terms_in_H},  $Y\in \cal H.$
Applying Holder's inequality for $p=1+\epsilon$ and $q=1 + \frac 1 {\epsilon}$ (satisfying $\frac 1 p + \frac 1 q =1$) yields that for any fixed $\epsilon >0$:
\begin{align}
     \Expect_{\pi\perp\tpi} & \Big[\Expect_\Q\big[L_\pi(X)L_{\tpi}(X) \Indc_{(X,\pi) \in \calE}\big] \Big] \nonumber \\
    &
    \, \leq \,  \Expect_{\pi\perp\tpi}
        \Bigg[
            \exp \Bigg( \big(1+o(1) \big) \, \rho^2 \sum_{1\leq k < \ell \leq m} \frac{ 1 + \Indc_{\{ d_{k\ell}(\sigma) < \frac n 2\}}} 2 \binom {d_{k\ell}(\sigma)} 2  + Y(\pi,\tpi) \Bigg)
        \Bigg] 
        \nonumber
    \\[0.5em]  
    & 
    \, \leq \, \Expect_{\pi\perp\tpi}
    \Bigg[
            \exp\Bigg( \big(1+o(1) \big) \, \rho^2(1+\epsilon) \sum_{1\leq k < \ell \leq m} \frac{ 1 + \Indc_{\{ d_{k\ell}(\sigma)< \frac n 2\}}} 2\binom {d_{k\ell}(\sigma)} 2 \Bigg)
    \Bigg]^{\frac 1 {1+\epsilon}} 
    \nonumber
    \\
    &  ~~~~~~~~~~~ \cdot  \Expect_{\pi\perp\tpi}
    \left[
        \exp\left(\left(1 + \frac 1 {\epsilon} \right)Y(\pi,\tpi)  \right)  
    \right]^{\frac{\epsilon}{1+\epsilon}}
    \nonumber
    \\[0.5em] 
    & 
    \, \leq \,  \Expect_{\pi\perp\tpi}
    \Bigg[
        \exp\Bigg( \big(1+o(1) \big) \, \rho^2(1+\epsilon) \sum_{1\leq k < \ell \leq m} \frac{ 1 + \Indc_{\{ d_{k\ell}(\sigma) < \frac n 2\}}} 2\binom {d_{k\ell}(\sigma)} 2  +o(1)\Bigg) 
    \Bigg]  
    \nonumber 
    \\[0.5em]
    &
    \, \leq \,  \Expect_{\pi\perp\tpi}
    \Bigg[
        \exp\Bigg(  \rho^2(1+2\epsilon) \sum_{1\leq k < \ell \leq m} \frac{ 1 + \Indc_{\{ d_{k\ell}(\sigma) < \frac n 2\}}} 2\binom {d_{k\ell}(\sigma)} 2  +o(1)\Bigg) 
    \Bigg]  \, .
     \label{eq:averaged}
\end{align}

\begin{remark}[Reduction to equality]\label{rem:rho-equality}
It is enough to prove~\Cref{thm:negative} when
\begin{equation}\label{eq:rho-equality}
        \rho^2 = \frac{8(1-\eta)}{m}\frac{\log n}{n} .
\end{equation}
Indeed, if $0\le \rho_0\le \rho_1$, then a sample from the planted model with correlation $\rho_1$ can be passed through the coordinatewise Gaussian channel
\[
        Y^k_e
        =
        \sqrt{\rho_0/\rho_1}\,X^k_e
        +
        \sqrt{1-\rho_0/\rho_1}\,Z^k_e,
        \qquad
        Z^k_e\stackrel{\mathrm{iid}}{\sim}N(0,1),
\]
with the same channel applied under the null. The null is mapped to the null, while the planted correlation is reduced from $\rho_1$ to $\rho_0$. Since total variation cannot increase under a Markov kernel, indistinguishability at the boundary correlation implies indistinguishability for all smaller nonnegative correlations.
\end{remark}

For the remainder of the proof of  ~\Cref{thm:negative} we assume that $\rho^2 = \frac{8(1-\eta)} m \frac {\log n} n$, 
and next show that the right-hand side of \eqref{eq:averaged} is $1 + o(1).$   Select $\epsilon$ such that  $(1-\eta)(1+2\epsilon) = 1 - {\eta}/2.$ 

Let $\etabar = \eta/2$ and $\rhobar^2 = \rho^2(1+2\epsilon)$ so that $\rhobar^2 = \frac{8(1-\etabar)} m \frac {\log n} n.$ Then \eqref{eq:averaged} yields 
\begin{align}   \label{eq:averaged_no_epsilon}
        \Expect_{\pi \perp \tpi}  &
        \Expect_\Q\big[L_\pi(X)L_{\tpi}(X) \Indc_{(X,\pi) \in \calE}\big] \notag \\[0.5em]
        & \leq  \Expect_{\pi\perp\tpi}
        \bigg[
            \exp\bigg(  \rhobar^2 \sum_{1\leq k < \ell \leq m}
            \frac{ 1 + \Indc_{\{ d_{k\ell}(\sigma) < \frac n 2\}}} 2\binom {d_{k\ell}(\sigma)} 2  +o(1) \bigg)
        \bigg]  \, \notag 
        \\[0.5em]
        & = \big(1+o(1)\big) \cdot \Expect_{\pi\perp\tpi}
        \bigg[
            \exp\bigg(  \rhobar^2 \sum_{1\leq k < \ell \leq m}
            \frac{ 1 + \Indc_{\{ d_{k\ell}(\sigma) < \frac n 2\}}} 2\binom {d_{k\ell}(\sigma)} 2  \bigg)
        \bigg] \, .
\end{align}
Now, let 
\[
    {\cal U} = \left\{ \sigma :  0 \leq d_{k\ell}(\sigma) \leq \frac{\sqrt{n}}{\log n} \mbox{ for } 1\leq k < \ell \leq m\right\} \, .
\]
Note that if $\sigma \in {\cal U}$, then 
\[
    \frac{1 + \Indc_{\{d_{k\ell}(\sigma) < \frac n 2 \}}} 2   \binom {d_{k\ell}(\sigma)} 2 
    \, \leq \, 
    \binom {d_{k\ell}(\sigma)} 2
    \, \leq \, 
    \frac 1 2 \frac n {\log^2 n} \, ,
\]
so that, using $\rhobar^2 \leq \frac 8 m \frac {\log n} n:$
\begin{align}
    \Expect_{\pi\perp\tpi} &
    \Bigg[
        \exp\Bigg( \rhobar^2 \!\! \sum_{1\leq k < \ell \leq m} \frac{ 1 + \Indc_{\{ d_{k\ell}(\sigma) < \frac n 2\}}} 2\binom {d_{k\ell}(\sigma)} 2\Bigg) \Indc_{\{\sigma\in\calU\}}
    \Bigg] \nonumber 
    \\[0.5em]
     & \leq  \exp
     \left(
        \binom m 2 \frac 1 2 \frac n {\log^2 n} \frac 8 m \frac {\log n} n 
    \right) \nonumber
    \\[0.5em]
     & =  \exp
    \left(
         \frac{2(m-1)}{\log n}   
    \right) \, .
     \label{eq:averaged_no_epsilon_U}
\end{align}

\begin{lemma}  \label{eq:f_bounds}
    For integers $d$ with $0\leq d \leq n$, let
    \[
        f(d) \coloneqq  - \log\big( d \, ! \big) +  \frac{2\big(1 + \Indc_{\{d < \frac n 2 \}} \big)  \cdot \binom d 2  \cdot (1- \etabar) \log n} n \, .
    \]
    Then, for all sufficiently large $n$, $f(d) \leq 0$ for $0\leq d \leq n$, and 
    \(
        \max \left\{ f(d) : \frac{\sqrt{n}}{\log n} \leq d \leq n\right\}
        \sim -\frac12\sqrt n \, .
    \)
\end{lemma}

\begin{proof}
Let 
\[ 
    g(d) \triangleq -\log(d!) +  \frac{4  \binom d 2 (1 - \etabar) \log n} n 
\]
so that $f(d) = g(d)$ for $1\leq d <  n / 2.$  The function $g$ satisfies 
\[ 
    g(0)=g(1)=0 
    ~~~~ \mbox{and} ~~~~
    g(d) - g(d-1) = -\log d + \frac{(d-1)4(1-\etabar)\log n} n, ~~ \mbox{ for }  d\geq 1 \, .
\]
Thus, for some $d^\star$, $g$ is decreasing over $0 \leq d \leq d^\star$ and increasing over $d\geq d^\star$. Since 
\[
    g\left(\frac n2\right) \sim - \frac{\etabar}{2} n\log n <0,
\]
it follows that $f(d)=g(d)\leq 0$ for $0 \leq d < n / 2.$

Next, by Stirling's formula, we have that
\(
    g\big(\frac{\sqrt n}{\log n}\big)
    \sim -\frac12\sqrt n
\)
and
\(
    g\big(\frac n2\big)
    \sim -\frac{\etabar}{2}n\log n .
\)
Since $g$ is decreasing and then increasing, the maximum of $g$ over 
$\sqrt n/\log n\le d<n/2$ is attained at one of the endpoints. Therefore
\[
    \max \Big\{ f(d) : \frac{\sqrt{n}}{\log n} \leq d < \frac n 2 \Big\} 
    =
    \max \Big\{ g\Big(\frac{\sqrt{n}}{\log n}\Big), \,  g\left( \frac n 2 \right) \! \Big\}
    \sim -\frac12\sqrt n ,
\]
where we used that $\etabar>0$ is fixed, so $\etabar n\log n \gg \sqrt n$.

Similarly, let 
\[
    h(d) \coloneqq -\log(d!) +  \frac{2  \binom d 2 (1 - \etabar) \log n} n 
\] 
so that $f(d) = h(d)$ for $\frac n 2 \leq d \leq  n.$ 
The function $h$ satisfies
\[
    h(d) - h(d-1) = -\log d + \frac{2(d-1)(1-\etabar)\log n} n 
    ~~ \mbox{ for }  d\geq 1 \, ,
\]
so, like $g$, it is decreasing for small $d$ and then increasing. Hence the maximum of $h$ over 
$n/2\le d\le n$ is attained at one of the endpoints. Moreover, Stirling's formula gives
\[
    h(n) \sim - \etabar \, n \log n,
    \qquad
    h\left(\frac n2\right)
    \sim -\frac{1+\etabar}{4} n\log n .
\]
Therefore, for some constant $c_{\etabar}>0$ depending only on $\etabar$,
\[
    \max \Big\{ f(d) :  \frac n 2 \leq d \leq n \Big\} 
    =
    \max \Big\{ h(d) : \frac n 2 \leq d \leq n \Big\}  
    \leq -c_{\etabar} n\log n .
\]
Combining the two ranges, and using again that $n\log n \gg \sqrt n$, we obtain as desired that
\[
    \max \Big\{ f(d): \frac{\sqrt n}{\log n} \leq d \leq n \Big\}
    \sim -\frac12\sqrt n  \, . \qedhere
\]
\end{proof}

Combining~\Cref{lemma:VM} and \Cref{eq:f_bounds} for $\rhobar^2 = \frac {8(1-\etabar)} m \frac {\log n} n$ yields:
\begin{align}
    \Expect_{\pi\perp\tpi} &
    \Bigg[
        \exp\Bigg( \rhobar^2 \sum_{1\leq k < \ell \leq m} \frac{ 1 + \Indc_{\{ d_{k\ell}(\sigma) < \frac n 2\}}} 2\binom {d_{k\ell}(\sigma)} 2\Bigg) \Indc_{\{\sigma\in\calU^c\}}
    \Bigg]
    \\[0.5em]
    &
    \leq
    \frac{1}{(n!)^{m-1}} 
    \sum_{(d_{k\ell}):\max_{k\ell}d_{k\ell}\geq\frac{\sqrt{n}}{\log n}} \!\!\!\!\!
     F((d_{k\ell})) 
     \exp
     \Bigg(  
        \Bigg[ \sum_{1\leq k < \ell \leq m}
            \frac{1 + \Indc_{\{d_{k\ell} < \frac n 2 \}}} 2   \binom {d_{k\ell}} 2
        \Bigg]  \, \rhobar^2
    \Bigg) 
    \nonumber
    \\[0.5em]
    & \stackrel{\mathrm{(a)}}{\leq}
    \sum_{(d_{k\ell}):\max_{k\ell}d_{k\ell}\geq\frac{\sqrt{n}}{\log n}}  \!\!\!\!
    \exp\Bigg(  
        \sum_{1\leq k <  \ell \leq m}  \Bigg[ - \frac 2 m \log(d_{k\ell}!) + 
      \frac{1 + \Indc_{\{d_{k\ell} < \frac n 2 \}}} 2   \binom {d_{k\ell}}{2} \,  \rhobar^2 \Bigg]  
    \Bigg)     
    \nonumber
    \\[0.5em]
    & =     \sum_{(d_{k\ell}):\max_{k\ell}d_{k\ell}\geq\frac{\sqrt{n}}{\log n}}  
    \exp\Bigg(  
        \frac 2 m  \sum_{1\leq k <  \ell \leq m} f(d_{k\ell})
    \Bigg)  
    \nonumber  
    \\[0.5em]  
    & \leq (n+1)^{\binom m 2} \exp\left(- \frac {2(1+o(1))} m  \times \frac {\sqrt{n}} 2  \right)  = o(1)  \, .  \label{eq:averaged_no_epsilon_Uc}
\end{align}
Here, (a) uses the bound in~\Cref{lemma:VM}.
Combining \eqref{eq:averaged_no_epsilon}, \eqref{eq:averaged_no_epsilon_U}, and \eqref{eq:averaged_no_epsilon_Uc} yields
\[
    \Expect_{\pi \perp \tpi} \ \Expect_\Q\big[L_\pi(X)L_{\tpi}(X) \Indc_{(X,\pi) \in \calE}\big] \leq 1 + o(1) \, .
\]
Returning to \eqref{eq:partial_progress} and using $\P(\calE)=1-o(1)$, we obtain
\[
        \Expect_\Q[(L')^2]\leq 1+o(1).
\]
Since $\Expect_\Q[L']=1$, this implies $\Expect_\Q[(L'-1)^2]=o(1)$ and hence $\TV(\P',\Q)=o(1)$.  Finally, by the high-probability estimate for $\calE$ and the standard bound between a measure and its conditioning on a high-probability event,
\[
        \TV(\P,\P')\leq \P(\calE^c)=e^{-\Omega(n)}.
\]
The triangle inequality gives $\TV(\P,\Q)=o(1)$, completing the proof of~\Cref{thm:negative}.

\appendix

\section{Supplementary proofs for achievability of detection}\label{app:achievability-tools}

\subsection{Proof of~\texorpdfstring{\Cref{lemma:HWGaussian}}{}}

Let $\lambda_1,\ldots,\lambda_d$ be the eigenvalues of $A$. By orthogonal invariance of the standard Gaussian distribution,
\(
    Z\stackrel{\mathrm{d.}}=\sum_{i=1}^d \lambda_i W_i^2,
\)
where $W_1,\cdots,W_d$ are independent standard Gaussians. Also
\[
    \Expect[Z] =\sum_{i=1}^d\lambda_i,
    ~~~~~~
    \norm{A}_{\mathrm F}^2=\sum_{i=1}^d\lambda_i^2,
    ~~~~~~
    \norm{A}_{\mathrm{op}}=\max_i |\lambda_i|.
\]
Fix $\theta>0$ such that $2\theta\norm{A}_{\mathrm{op}}<1$. Then
\begin{align*}
        \Expect\exp\{\theta(Z-\Expect Z)\}
        &=
        \prod_{i=1}^d
        \exp\{-\theta\lambda_i\}(1-2\theta\lambda_i)^{-1/2} .
\end{align*}
Using the power series for $\log(1+x)$ and the inequality $|2\theta\lambda_i|<1$,
\[
        -\frac12\log(1-2\theta\lambda_i)-\theta\lambda_i
        \le
        \frac{\theta^2\lambda_i^2}{1-2\theta\norm{A}_{\mathrm{op}}}.
\]
Therefore
\[
        \Expect\exp\{\theta(Z-\Expect Z)\}
        \le
        \exp\left(
        \frac{\theta^2\norm{A}_{\mathrm F}^2}
        {1-2\theta\norm{A}_{\mathrm{op}}}
        \right).
\]
By Chernoff's bound,
\[
        \Prob(Z-\Expect Z\ge t)
        \le
        \exp\left(
        \frac{\theta^2\norm{A}_{\mathrm F}^2}
        {1-2\theta\norm{A}_{\mathrm{op}}}
        -\theta t
        \right).
\]
Choosing
\(
        \theta
        =
        \frac{t}{2\left(\norm{A}_{\mathrm F}^2+\norm{A}_{\mathrm{op}}t\right)}
\)
gives \eqref{eq:HWGaussian-one-sided}. Finally, \eqref{eq:HWGaussian-min-form} follows from
\[
        \frac{1}{a+b}
        \ge
        \min\left\{\frac{\gamma}{a},\frac{1-\gamma}{b}\right\},
        \qquad a,b>0.
\]

\subsection{Proof of~\texorpdfstring{\Cref{lem:edge-quadratic-eigs}}{}}

Let $\mathbf 1$ denote the all ones vector in $\mathbb{R}^m$ and define
\(
    H\coloneqq \frac12(\mathbf 1\mathbf 1^\top-I_m).
\)
Then $Z_\rho=Y^\top H Y$.  We use the representation
\[
        Y\stackrel{\mathrm d}=
        \begin{pmatrix}
        \sqrt{\rho}\,\mathbf 1 & \sqrt{1-\rho}\,I_m
        \end{pmatrix}W ,
        \qquad
        W\sim\mathcal N(0,I_{m+1}).
\]
Indeed, the covariance of the right-hand side is
$\rho\mathbf 1\mathbf 1^\top+(1-\rho)I_m=\Sigma_\rho$.
Writing
\[
    B\coloneqq
    \begin{pmatrix}
    \sqrt{\rho}\,\mathbf 1 & \sqrt{1-\rho}\,I_m
    \end{pmatrix},
    ~~~~~~
    \bar A_\rho\coloneqq B^\top H B,
\]
we have
\(
    Z_\rho\stackrel{\mathrm d.} = W^\top \bar A_\rho W .
\)
It remains to identify the nonzero eigenvalues of $\bar A_\rho$.
First, set
\[
    v_+\coloneqq
    \begin{pmatrix}
    m\sqrt{\rho}\\
    \sqrt{1-\rho}\,\mathbf 1
    \end{pmatrix}.
\]
Since $Bv_+=(1+(m-1)\rho)\mathbf 1$, and $H\mathbf 1=\frac{m-1}{2}\mathbf 1$, and
$B^\top\mathbf 1=v_+$, it follows that
\[
    \bar A_\rho v_+
    =
    \frac{m-1}{2}\big(1+(m-1)\rho\big)v_+ .
\]
Thus,
\(
    \lambda_+(\rho)
    =
    \frac{m-1}{2}\big(1+(m-1)\rho\big)
\)
is an eigenvalue of $\bar A_\rho$.
Next, if $u\perp\mathbf 1$ and
\(
    v_u\coloneqq
    \begin{pmatrix}
        0 \\
        u
    \end{pmatrix},
\)
then $Bv_u=\sqrt{1-\rho}\,u$ and $Hu=- u/2$.  Therefore
\[
    \bar A_\rho v_u
    =
    B^\top H Bv_u
    =
    -\frac{1-\rho}{2}v_u .
\]
The space of such $u$ has dimension $m-1$, so
\(
    \lambda_-(\rho)=-\frac{1-\rho}{2}
\)
is an eigenvalue with multiplicity $m-1$.
Finally,
\[
        v_0\coloneqq
        \begin{pmatrix}
        -\sqrt{1-\rho}\\
        \sqrt{\rho}\,\mathbf 1
        \end{pmatrix}
\]
satisfies $Bv_0=0$, and hence $\bar A_\rho v_0=0$.  The vector
$v_0$, the vector $v_+$, and the subspace
$\{(0,u)^\top:u\perp\mathbf 1\}$ account for all $m+1$
orthogonal directions.  Thus the only nonzero eigenvalues of
$\bar A_\rho$ are $\lambda_+(\rho)$ with multiplicity one and
$\lambda_-(\rho)$ with multiplicity $m-1$.

Since $\bar A_\rho$ is symmetric, an orthogonal change of coordinates
diagonalizes the quadratic form, and this change of coordinates preserves
the standard Gaussian law of $W$.  Deleting the coordinate corresponding
to the zero eigenvalue gives
\(
    Z_\rho\stackrel{\mathrm d}=g^\top A_\rho g,
\)
for
\(
    g\sim\mathcal N(0,I_m),
\)
and $A_\rho$ has the claimed eigenvalues and multiplicities.

\section{Supplementary proofs for impossibility of detection} \label{app:matrix-properties}

\subsection{Proof of~\texorpdfstring{\Cref{lemma:K_properties}}{}}

(a) The equality $K'_{\pi}=K'_{\tpi}$ follows from the fact that if $e$ is an edge with $e\in \binom {O_{k\ell}} 2$ then
$e=\{i,j\}$ such that $\pi_{k\ell}(i)=\tpi_{k\ell}(i)$ and $\pi_{k\ell}(j)=\tpi_{k\ell}(j)$ so $\pi_{k\ell}(e)=\tpi_{k\ell}(e).$

\medskip 

(b) The first inequality follows directly from the definition and the second from part (a).

\medskip 

(c) 
Using the definition \eqref{eq:K-entry},
\[
(K_\pi)_{(k,e),(\ell,f)} = \ind{k\ne \ell}\,\ind{\pi_{k\ell}(e)=f} 
 = \ind{\ell\ne k}\,\ind{\pi_{\ell k}(f)=e}  = (K_\pi)_{(\ell,f),(k,e)},
\]
so $K_{\pi}$ is symmetric.  Similarly,
\begin{align*}
    (K'_\pi)_{(k,e),(\ell,f)} &= \ind{k\ne \ell}\,\ind{\pi_{k\ell}(e)=f}\ind{e\in \binom{O_{k\ell}}{2}}  \\
    & = \ind{\ell\ne k}\,\ind{\pi_{\ell k}(f)=e}\ind{f \in \binom{O_{\ell k}}{2}}  = (K'_\pi)_{(\ell,f),(k,e)}
\end{align*}
so $K_{\pi}'$ is also symmetric.

\medskip

(d) Observe that
\begin{align*}
    (K_\pi K_{\tpi})_{(k,e),(k,e)}
    &=
    \sum_{\ell=1}^m\sum_{f\in\En}
    (K_\pi)_{(k,e),(\ell,f)}
    (K_{\tpi})_{(\ell,f),(k,e)} \\
    &=
    \sum_{\ell\ne k}\sum_{f\in\En}
    \ind{\pi_{k\ell}(e)=f}\,
    \ind{\tpi_{\ell k}(f)=e} \\
    &=
    \sum_{\ell\ne k}
    \ind{\pi_{k\ell}(e)=\tpi_{k\ell}(e)}.
\end{align*}
Taking the trace gives
\begin{align*}
    \Tr(K_\pi K_{\tpi})
    &=
    \sum_{k=1}^m\sum_{e\in\En}
    (K_\pi K_{\tpi})_{(k,e),(k,e)} \\
    &=
    \sum_{k=1}^m\sum_{\ell\ne k}
    \sum_{e\in\En}
    \ind{\pi_{k\ell}(e)=\tpi_{k\ell}(e)} \\
    &=
    \sum_{k=1}^m\sum_{\ell\ne k}D_{k\ell}(\pi,\tpi).
\end{align*}
Finally, $D_{k\ell}=D_{\ell k}$ and $D_{k\ell}(\pi,\tpi)=D_{k\ell}(\sigma)$, so the last sum equals $2\sum_{1\le k<\ell\le m}D_{k\ell}(\sigma)$.

\medskip 

(e)  Observe that
\begin{align*}
    (K'_\pi K_{\tpi})_{(k,e),(k,e)}
    &=
    \sum_{\ell=1}^m\sum_{f\in\En}
    (K'_\pi)_{(k,e),(\ell,f)}
    (K_{\tpi})_{(\ell,f),(k,e)} \\
    &=
    \sum_{\ell\ne k}\sum_{f\in\En}
    \ind{\pi_{k\ell}(e)=f}\,\ind{e\in \binom{O_{k\ell}}{2}}
    \ind{\tpi_{\ell k}(f)=e} \\
    &=
    \sum_{\ell\ne k} \ind{e\in \binom{O_{k\ell}}{2}} 
    \ind{\pi_{k\ell}(e)=\tpi_{k\ell}(e)}.
\end{align*}
Taking the trace gives
\begin{align*}
    \Tr(K'_\pi K_{\tpi})
    &=
    \sum_{k=1}^m\sum_{e\in\En}
    (K'_\pi K_{\tpi})_{(k,e),(k,e)} \\
    &=
    \sum_{k=1}^m\sum_{\ell\ne k}
    ~~\sum_{e\in\En\cap \binom{O_{k\ell}}{2}}
    \ind{\pi_{k\ell}(e)=\tpi_{k\ell}(e)} \\
        &=
    \sum_{k=1}^m\sum_{\ell\ne k}
    \Big|\En\cap \binom{O_{k\ell}}{2} \Big| =
    \sum_{k=1}^m\sum_{\ell\ne k}  \binom {d_{k\ell}(\sigma)} 2  \Indc_{\{ d_{k\ell}(\sigma)\geq \frac n 2\}}.
\end{align*}
Finally, $d_{k\ell}(\sigma)=d_{\ell k}(\sigma)$, so the last sum equals $2\sum_{1\le k<\ell\le m}\binom {d_{k\ell}(\sigma)} 2  \Indc_{\{ d_{k\ell}(\sigma)\geq \frac n 2\}}.$  Proof of the other equalities is similar (and also $\Tr(AB)=\Tr(BA)$ for matrices of appropriate dimensions in general).

\subsection{Proof of~\texorpdfstring{\Cref{lem:trace-expansions}}{}}

For part~\textup{(a)}, expanding the diagonal entries of $K^2_\pi K_{\tpi}$ gives

\begin{align*}
    \big(K^2_{\pi} \, K_{\tpi}\big)_{(k,e),(k,e)} 
    \, =
    \sum_{\ell,k' \in [m]}
    \sum_{~f,e',f'\in\En}
    (K_\pi)_{(k,e),(\ell,f)} \, 
    (K_{\pi})_{(\ell,f),(k',e')} \, 
    (K_{\tpi})_{(k',e'),(k,e)} .
\end{align*}
Using \eqref{eq:K-entry}, this becomes
\begin{align*}
    \big(K^2_\pi \, K_{\tpi}\big)_{(k,e),(k,e)}  
    \, =
    \sum_{\ell,k' \in [m]}
    \sum_{~f,e'\in\En}
    \ind{k\ne\ell,\ \ell\ne k',\ k'\ne k} \cdot 
    \ind{\pi_{k\ell}(e)=f} \cdot 
    \ind{\pi_{\ell k'}(f)=e'} \cdot 
        \ind{\tpi_{k' k}(e')=e}\,
    .
\end{align*}
For fixed $k,\ell,k'$ and $e$, the first two edge indicators uniquely determine
$f=\pi_{k\ell}(e)$ and $e'=\tpi_{\ell k'}(f).$ 
Hence the final condition is precisely
\[
        \big(
        \tpi_{k' k}
        \circ \pi_{\ell k'}
        \circ \pi_{k \ell}
        \big)(e)=e
\]
for $\ell,k'$ such that $k,\ell,k'$ are distinct.
Summing over all diagonal coordinates $(k,e)$ gives part~\textup{(a)}.

Part~\textup{(b)} is similar, expanding the diagonal of $K_{\pi} K_{\tpi} K_{\pi} K_{\tpi}$; the first three indicators determine the intermediate edge variables, and the final indicator becomes
$\big(\tpi_{\ell'k}\circ\pi_{k'\ell'}\circ\tpi_{\ell k'}\circ\pi_{k\ell}\big)(e)=e$
where the sum is over $k,\ell,k',\ell'$ such that $\{k,k'\}\cap \{\ell, \ell'\}=\emptyset.$

Part~\textup{(c)} is identical to part (b).  Expanding the diagonal of $K_\pi K_\pi K_\pi K_{\tpi}$; the final indicator becomes
$\big(\tpi_{\ell'k}\circ\pi_{k'\ell'}\circ\pi_{\ell k'}\circ\pi_{k\ell}\big)(e)=e$.

\subsection{Proof of~\texorpdfstring{\Cref{lem:exp-Cj-square}}{}} \label{sec-proof-of-lem:exp-Cj-square}

A classical bound on the tail of the distribution of $C_j$ is \cite{arratia2003logarithmic}
\[
    \mathbb{P}(C_j \ge k) \le \frac{1}{j^k k!} \, .
\]
(It follows from $\mathbb{P}(C_j \geq k) = \mathbb{P}\big(\binom {C_j} k \geq 1\big)$, Markov's inequality, and the fact $\Expect\big[  \binom {C_j} k  \big] \leq \frac 1 {j^k k!}$ with equality for $jk\leq n$ \cite{arratia2003logarithmic,Feller1957}).
Also $C_j \leq n/j.$  
We shall use the following summation by parts identity: if \(Z\) is a nonnegative integer-valued random variable and \(a_0,a_1,\ldots\) is increasing, then
\[
    \mathbb E[a_Z]
    =
    a_0+\sum_{k\ge1}(a_k-a_{k-1})\, \mathbb P(Z\ge k).
\]
For part \(\mathrm{(a)}\), first assume \(\theta_n\ge0\). Then
\[
\begin{aligned}
    \mathbb E[e^{\theta_n C_j}]-1
    &\leq
    \sum_{k\ge1}
    \left(e^{\theta_n k}-e^{\theta_n(k-1)}\right)
    \frac1{k!}  \\
    &=
    (e^{\theta_n}-1)
    \sum_{k\ge1}
    \frac{e^{\theta_n(k-1)}}{k!}
    =
    O(\theta_n)=o(1).
\end{aligned}
\]
If \(\theta_n<0\), then
\[
    0\le 1-\mathbb E[e^{\theta_n C_j}]
    \le \mathbb E[e^{|\theta_n|C_j}]-1=o(1),
\]
so part \(\mathrm{(a)}\) follows.

To complete the proof of the lemma, the same reduction lets us consider only \(\theta_n\ge0\). Since \(\theta_n=o(\log n/n)\), by monotonicity it suffices to take $\theta_n =  {\log n}/{(2n)}$
and show that $\mathbb{E}[e^{\theta_n C_j^2}] = 1 + o(1)$.
Using the same summation identity,
\[
    \mathbb{E}[e^{\theta_n C_j^2}]-1
    \le
    \sum_{k=1}^{n}
    \left(e^{\theta_n k^2}-e^{\theta_n(k-1)^2}\right)\frac1{k!}.
\]
Let  $M_n \coloneqq \left\lfloor \frac{n}{(\log n)^2} \right\rfloor$ and split the last sum into
\[
    S_{n,1} \coloneqq \sum_{k=1}^{M_n}
    \left(e^{\theta_n k^2}-e^{\theta_n(k-1)^2}\right)\frac1{k!},
    ~~~~~~
    S_{n,2} \coloneqq \sum_{k=M_n+1}^{n}
    \left(e^{\theta_n k^2}-e^{\theta_n(k-1)^2}\right)\frac1{k!}.
\]

\textit{Small $k$ part.}
For $k\le M_n$ we have $k^2\le M_n k$ and for $x\geq 0$ we have $1-e^{-x} \leq x$, so that
\[
    e^{\theta_n k^2}-e^{\theta_n(k-1)^2} = (1-e^{-\theta_n(2k-1)})e^{\theta_nM_nk}
    \le
    2\theta_n k\,e^{\theta_n M_n k}.
\]
Therefore
\[
    S_{n,1}
    \le
    2\theta_n
    \sum_{k\ge1}
    \frac{k e^{\theta_n M_n k}}{k!}
    =
    O(\theta_n)=o(1).
\]

\textit{Large $k$ part.}
For \(S_{n,2}\), we use
\(e^{\theta_n k^2}-e^{\theta_n(k-1)^2}\le e^{\theta_n k^2}\).
Using $k!>(k/e)^k$,
\[
    \frac{e^{\theta_n k^2}}{k!}
    < \exp \left(\theta_n k^2 + k - k\log k\right)
    \eqqcolon  e^{f_n(k)}.
\]
A direct computation shows:
\begin{itemize}
\item $f_n''(k)=2\theta_n - 1/k$ has exactly one zero at
      $k_0 = n/\log n$;
\item $M_n < k_0 < n$ for large $n$;
\item $f_n'(M_n)<0$, $f_n'(k_0)<0$, and $f_n'(n)=0$.
\end{itemize}
Thus $f_n'$ decreases on $[M_n,k_0]$ and increases on $[k_0,n]$, remaining $\le 0$ throughout. Hence $f_n$ is strictly decreasing on $[M_n,n]$, so
\[
    f_n(k) \le f_n(M_n)\qquad\text{for all }k\in[M_n,n] \, .
\]
A direct asymptotic expansion gives
\[
    f_n(M_n)
    = -\frac{n}{\log n} + O\!\left(\frac{n\log\log n}{(\log n)^2}\right)
    \le -\frac{1}{2}\,\frac{n}{\log n}
\]
for all large $n$. Therefore
\[
    S_{n,2}
    \le
    n \exp\!\left(-\frac{1}{2}\,\frac{n}{\log n}\right) 
    \to 0.
\]
Combining the two estimates gives
\[
\mathbb{E}[e^{\theta_n C_j^2}]-1\le S_{n,1}+S_{n,2}=o(1),
\]
which proves part \(\mathrm{(b)}\).

\section{On the test statistic}\label{app:GLRT-derivation}

In this appendix we show that the statistic $T$ defined in~\eqref{eq:T-def} is approximately a generalized likelihood ratio statistic. 
By definition, the generalized likelihood ratio is $\max_{\pi} L_\pi(X)$ and thresholding any (strictly) increasing function of the generalized likelihood ratio is a generalized likelihood ratio test (GLRT).
By \eqref{eq:LLR},  $L_\pi(X)$ is a monotone increasing function of 
$X^\top(I_{mN} - \Sigma_\pi^{-1})X$ where $\Sigma_{\pi} = I_{mN} + \rho K_{\pi}.$
Calculating as in the proof of  Lemma \ref{lemma:K_properties} it can be shown that
\begin{align*}
K_\pi^2 = (m-1)I_{mN} +  (m-2)K_{\pi}
\end{align*}
from which it can be readily checked that
\begin{align*}
\Sigma_\pi^{-1} = \frac{1+(m-2)\rho}{(1-\rho)(1+(m-1)\rho)}I_{mN}
- \frac{\rho}{(1-\rho)(1+(m-1)\rho)}K_{\pi}.
\end{align*}
(This may look complicated but it can be derived by considering the case $n=1$ because up to a permutation of rows and columns depending on $\pi$, the matrix $\Sigma_{\pi}$ is block diagonal with $\binom n 2$ blocks of the form $\Sigma$ defined in \eqref{eq:small_sigma}.)
Therefore, using the fact  $ X^\top K_{\pi}X = 2T_\pi(X),$ we get
\begin{align*}
    X^\top(I_{mN} - \Sigma_\pi^{-1})X =   \frac{ 2\rho T_{\pi}(X)  -(m-1)\rho^2  \|X\|_2^2  }{(1-\rho)(1+(m-1)\rho)}.
\end{align*}
Since the coefficient of $\|X\|_2^2$ includes $\rho^2$ which is much smaller than the coefficient $2\rho$ of $T_{\pi}(X)$, we conclude that $L_{\pi}(X)$ is approximately a monotone increasing function of $T_{\pi}(X).$  Thus, $T$ is approximately a generalized likelihood ratio statistic. (In addition, selecting $\pi$ to maximize $T_{\pi}(X)$ is the maximum likelihood estimator of the true alignment in case it is known that $H_1$ is the true hypothesis.)

\subsection*{Acknowledgments}
This work was supported by NSF under Grant CCF 19-00636. The authors acknowledge the use of AI-based tools for assistance with literature search, editorial refinement, and improving the clarity of exposition.

\bibliographystyle{alpha}
\bibliography{bibliography}

\end{document}